\documentclass{amsart}
\usepackage{amssymb, amsmath, mathrsfs, verbatim, url, bbm}

\theoremstyle{plain}
\newtheorem{theorem}{Theorem}[section]
\newtheorem{lemma}[theorem]{Lemma}
\newtheorem{proposition}[theorem]{Proposition}
\newtheorem{corollary}[theorem]{Corollary}

\newtheorem{observation}[theorem]{Observation}

\theoremstyle{definition}
\newtheorem{definition}[theorem]{Definition}
\newtheorem{example}[theorem]{Example}

\theoremstyle{remark}
\newtheorem*{remark}{Remark}

\newtheorem{question}[theorem]{Question}

\numberwithin{equation}{section}

\DeclareMathOperator{\upset}{\uparrow\!}

\DeclareMathOperator{\dom}{dom}
\DeclareMathOperator{\ran}{ran}
\DeclareMathOperator{\cf}{cf}

\DeclareMathOperator{\Clop}{Clop}
\DeclareMathOperator{\diam}{diam}
\DeclareMathOperator{\Lim}{Lim}
\newcommand{\Aut}[1]{\mathrm{Aut}(#1)}

\newcommand{\nbd}{\nobreakdash}
\newcommand{\homeo}{\cong}
\newcommand{\powset}[1]{\mathcal{P}(#1)}
\newcommand{\card}[1]{\lvert #1\rvert}

\newcommand{\la}{\langle}
\newcommand{\ra}{\rangle}
\newcommand{\ord}{{\mathrm On}}
\newcommand{\cardc}{\mathfrak{c}}
\newcommand{\cardp}{\mathfrak{p}}
\newcommand{\inv}[1]{#1^{-1}}
\newcommand{\closure}[1]{\overline{#1}}
\newcommand{\restrict}{\upharpoonright}

\newcommand{\cell}[1]{c(#1)}
\newcommand{\density}[1]{d(#1)}
\newcommand{\weight}[1]{w(#1)}
\newcommand{\ow}[1]{Nt(#1)}

\newcommand{\opi}[1]{\pi Nt(#1)}
\newcommand{\ochar}[1]{\chi Nt(#1)}

\newcommand{\okchar}[1]{\chi_K Nt(#1)}

\newcommand{\piweight}[1]{\pi(#1)}
\newcommand{\character}[1]{\chi(#1)}

\newcommand{\reals}{\mathbb{R}}

\newcommand{\integers}{\mathbb{Z}}
\newcommand{\Fn}[2]{\mathrm{Fn}(#1,\,#2)}
\newcommand{\wma}{we may assume}
\newcommand{\Wma}{We may assume}
\newcommand{\op}{\mathrm{op}}

\newcommand{\omegaop}{$\omega^\op$\nbd-like}
\newcommand{\kappaop}{$\kappa^\op$\nbd-like}
\newcommand{\cardop}[1]{$#1^\op$\nbd-like}
\newcommand{\elemsub}{\prec}

\newcommand{\supp}[1]{\mathrm{supp}(#1)}

\newcommand{\pisw}[1]{\pi sw(#1)}
\newcommand{\tightness}[1]{t(#1)}

\newcommand{\mcA}{\mathcal{ A}}
\newcommand{\mcB}{\mathcal{ B}}
\newcommand{\mcC}{\mathcal{ C}}
\newcommand{\mcD}{\mathcal{ D}}
\newcommand{\mcE}{\mathcal{ E}}
\newcommand{\mcF}{\mathcal{ F}}
\newcommand{\mcG}{\mathcal{ G}}
\newcommand{\mcH}{\mathcal{ H}}
\newcommand{\mcI}{\mathcal{ I}}

\newcommand{\mcL}{\mathcal{ L}}

\newcommand{\mcN}{\mathcal{ N}}

\newcommand{\mcQ}{\mathcal{ Q}}
\newcommand{\mcR}{\mathcal{ R}}
\newcommand{\mcS}{\mathcal{ S}}
\newcommand{\mcT}{\mathcal{ T}}
\newcommand{\mcU}{\mathcal{ U}}
\newcommand{\mcV}{\mathcal{ V}}
\newcommand{\mcW}{\mathcal{ W}}

\newcommand{\mbP}{\mathbb{ P}}
\newcommand{\mbQ}{\mathbb{ Q}}

\newcommand{\mbT}{\mathbb{ T}}
\newcommand{\msS}{\mathscr{S}}

\begin{document}

\title{Noetherian types of homogeneous compacta and dyadic compacta}

\author{David Milovich}

\date{May 15, 2007}

\address{University of Wisconsin-Madison Mathematics Dept.}
\email{milovich@math.wisc.edu}
\thanks{Support provided by an NSF graduate fellowship.}

\begin{abstract}
The Noetherian type of a space is the least $\kappa$ such that it has a base that is $\kappa$\nbd-like with respect to containment.  Just as all known homogeneous compacta have cellularity at most $\cardc$, they satisfy similar upper bounds in terms of Noetherian type and related cardinal functions.  We prove these and many other results about these cardinal functions.  For example, every homogeneous dyadic compactum has Noetherian type $\omega$.  Assuming GCH, every point in a homogeneous compactum $X$ has a local base that is $\cell{X}$\nbd-like with respect to containment.  If every point in a compactum has a well-quasiordered local base, then some point has a countable local $\pi$\nbd-base.
\end{abstract}

\maketitle

\section{Introduction}

Van Douwen's Problem (see Kunen~\cite{kunenlhc}) asks whether there is a homogeneous compactum of cellularity exceeding $\cardc$.  (See Engelking~\cite{engelking}, Juh\'asz~\cite{juhasz}, and Kunen~\cite{kunenst} for all undefined terms.  In particular, recall that $\weight{\cdot}$, $\piweight{\cdot}$, $\character{\cdot}$, $\pi\character{\cdot}$, $\density{\cdot}$, $\cell{\cdot}$, and $\tightness{\cdot}$ respectively denote weight, $\pi$\nbd-weight, character, $\pi$\nbd-character, density, cellularity, and tightness of topological spaces.)
A homogeneous compactum of cellularity $\cardc$ exists by Maurice~\cite{maurice}, but van Douwen's Problem remains open in all models of ZFC.

\begin{definition}
We say that a homogeneous compactum is \emph{exceptional} if it is not homeomorphic to a product of dyadic compacta and first countable compacta.
\end{definition}

By Arhangel${}^\prime$ski\u\i's Theorem, first countable spaces have size at most $\cardc$; dyadic compacta are ccc.  Since the cellularity of a product space equals the supremum of the cellularities of its finite subproducts (see p. 107 of \cite{juhasz}), all nonexceptional homogeneous compacta have cellularity at most $\cardc$.   To the best of the author's knowledge, there are only two classes of examples of exceptional homogeneous compacta (see \cite{milovicha}); these two kinds of spaces have cellularities $\omega$ and $\cardc$.

We investigate several cardinal functions defined in terms order\nbd-theoretic base properties.  Just like cellularity, these functions have upper bounds when restricted to the class of known homogeneous compacta.  Moreover, GCH implies that one of these functions is a lower bound on cellularity when restricted to homogeneous compacta.

\begin{definition}
Given a cardinal $\kappa$, define a poset to be $\kappa$\nbd-\emph{like} ($\kappa^{\op}$\nbd-\emph{like}) if no element is above (below) $\kappa$\nbd-many elements.  Define a poset to be \emph{almost} $\kappa^{\op}$\nbd-\emph{like} if it has a \kappaop\ dense subset.
\end{definition}

In the context of families of subsets of a topological space, we will always implicitly order by inclusion.

\begin{definition}
Given a space $X$. let the \emph{Noetherian type} of $X$, or $\ow{X}$, be the least $\kappa\geq\omega$ such that $X$ has a base that is \kappaop.  Analogously define \emph{Noetherian $\pi$\nbd-type} in terms of $\pi$\nbd-bases and denote it by $\pi Nt$.  Given a subset $E$ of $X$, let the \emph{local Noetherian type} of $E$ in $X$, or $\ochar{E,X}$, be the least $\kappa\geq\omega$ such that there is a \kappaop\ neighborhood base of $E$.  Given $p\in X$, let the local Noetherian type of $p$, or $\ochar{p,X}$, be $\ochar{\{p\},X}$.  Let the local Noetherian type of $X$, or $\ochar{X}$, be the supremum of the local Noetherian types of its points.  Let the \emph{compact Noetherian type} of $X$, or $\okchar{X}$, be the supremum of the local Noetherian types of its compact subsets.  We call $Nt$, $\pi Nt$, $\chi Nt$, and $\chi_K Nt$ \emph{Noetherian cardinal functions}.
\end{definition}

Noetherian type and Noetherian $\pi$\nbd-type were introduced by Peregudov~\cite{peregudov97}.  Preceding this introduction are several papers by Peregudov, Shapirovski\u\i\, and Malykhin~\cite{malykhin,peregudov76c,peregudov76p,peregudovS76} about $\min\{\ow{\cdot},\omega_2\}$ and $\min\{\opi{\cdot},\omega_2\}$ (using different terminologies).  Also, Dow and Zhou~\cite{dow} showed that $\beta\omega\setminus\omega$ has a point with local Noetherian type $\omega$.

\begin{observation}\label{OBShomogbound}
Every known homogeneous compactum $X$ satisfies the following.
\begin{enumerate}
\item\label{enumobsow} $\ow{X}\leq\cardc^+$.
\item\label{enumosbopi} $\opi{X}\leq\omega_1$.
\item\label{enumobsochar} $\ochar{X}=\omega$.
\item\label{enumobsokchar} $\okchar{X}\leq\cardc$.
\end{enumerate}
\end{observation}

We justify this observation in Section~\ref{SECupper}, except that we postpone the case of homogeneous dyadic compacta to Section~\ref{SECdyadic}, where we investigate Noetherian cardinal functions on dyadic compacta in general.  The results relevant to Observation~\ref{OBShomogbound} are summarized by the following theorem.

\begin{theorem}\label{THMdyadicmain}
Suppose $X$ is a dyadic compactum.  Then $\opi{X}=\okchar{X}=\omega$.  Moreover, if $X$ is homogeneous, then $\ow{X}=\omega$.
\end{theorem}

Also in Section~\ref{SECdyadic}, we generalize the above theorem to continuous images of products of compacta with bounded weight; we also prove the following.

\begin{theorem}\label{THMdyadicspectrum}
The class of Noetherian types of dyadic compacta includes $\omega$, excludes $\omega_1$, includes all singular cardinals, and includes $\kappa^+$ for all cardinals $\kappa$ with uncountable cofinality.
\end{theorem}

Section~\ref{SECreflectcone} investigates to what extent a technical property of free boolean algebras that is crucial to Section~\ref{SECdyadic} holds in other boolean algebras.  In Section~\ref{SEChomogloc}, we prove several results about the local Noetherian types of all homogeneous compacta, known and unknown, including the following theorem.

\begin{theorem}[GCH]\label{THMgchocharcell}
If $X$ is a homogeneous compactum, then $\ochar{X}\leq\cell{X}$.
\end{theorem}

\section{Observed upper bounds on Noetherian cardinal functions}\label{SECupper}
First, we note some very basic facts about Noetherian cardinal functions.

\begin{definition}
Given a subset $E$ of a product $\prod_{i\in I}X_i$ and $\sigma\in[I]^{<\omega}$, we say that $E$ has \emph{support} $\sigma$, or $\supp{E}=\sigma$, if $E=\inv{\pi_\sigma}\pi_\sigma[E]$ and $E\not=\inv{\pi_\tau}\pi_\tau[E]$ for all $\tau\subsetneq\sigma$.
\end{definition}

\begin{theorem}\label{THMowprodsup}
Given a point $p$ and a compact subset $K$ of a product space $X=\prod_{i\in I}X_i$, we have the following relations.
\begin{align*}
\ow{X} &\leq\sup_{i\in I}\ow{X_i}\\
\opi{X} &\leq\sup_{i\in I}\opi{X_i}\\
\ochar{p,X} &\leq\sup_{i\in I}\ochar{p(i),X_i}\\
\ochar{K,X} &\leq\sup_{\sigma\in[I]^{<\omega}}\okchar{\pi_\sigma[K],\pi_\sigma[X]}
\end{align*}
\end{theorem}
\begin{proof}
See Peregudov~\cite{peregudov97} for a proof of the first relation. That proof can be easily modified to demonstrate the next two relations.  Let us prove the last relation.  For each $\sigma\in[I]^{<\omega}$, set $\kappa_\sigma=\ochar{\pi_\sigma[K],\pi_\sigma[X]}$ and let $\mcA_\sigma$ be a \cardop{\kappa_\sigma}\ neighborhood base of $\pi_\sigma[K]$.  For each $\sigma\in[I]^{<\omega}$, let $\mcB_\sigma$ denote the set of sets of the form $\inv{\pi_\sigma}U$ where $U\in\mcA_\sigma$ and $\supp{U}=\sigma$.  Note that if $U\in\mcA_\sigma$ and $\supp{U}\subsetneq\sigma$, then there exists $\tau\subsetneq\sigma$ and $V\in\mcA_\tau$ such that $\inv{\pi_\tau}V\subseteq\inv{\pi_\sigma}U$.  Moreover, for any minimal such $\tau$, we have $\inv{\pi_\tau}V\in\mcB_\tau$.

Set $\mcB=\bigcup_{\sigma\in[I]^{<\omega}}\mcB_\sigma$.  By compactness, $\mcB$ is a neighborhood base of $K$.  Moreover, if $\sigma,\tau\in[I]^{<\omega}$ and $\mcB_\sigma\ni U\subseteq V\in\mcB_\tau$, then $\sigma=\supp{U}\supseteq\supp{V}=\tau$; hence, given $U$, there are at most $(\sup_{\tau\subseteq\sigma}\kappa_\tau)$\nbd-many possibilities for $V$.  Thus, $\mcB$ is \cardop{(\sup_{\sigma\in[I]^{<\omega}}\kappa_\sigma)}\ as desired.
\end{proof}

\begin{question}
Do there exist spaces $X$ and $Y$ such that $\okchar{X\times Y}$ exceeds $\okchar{X}\okchar{Y}$?
\end{question}

\begin{lemma}\label{LEMkappasizeposet}
Every poset $P$ is almost $\card{P}^{\op}$\nbd-like.
\end{lemma}
\begin{proof}
Let $\kappa=\card{P}$ and let $\la p_\alpha\ra_{\alpha<\kappa}$ enumerate $P$.  Define a partial map $f\colon\kappa\to P$ as follows.  Suppose $\alpha<\kappa$ and we have a partial map $f_\alpha\colon\alpha\to P$.  If $\ran f_\alpha$ is dense in $P$, then set $f_{\alpha+1}=f_\alpha$.  Otherwise, set $\beta=\min\{\delta<\kappa:p_\delta\not\geq q\text{ for all }q\in\ran f_\alpha\}$ and let $f_{\alpha+1}$ be the smallest map extending $f_\alpha$ such that $f_{\alpha+1}(\alpha)=p_\beta$.  For limit ordinals $\gamma\leq\kappa$, set $f_\gamma=\bigcup_{\alpha<\gamma}f_\alpha$.  Then $f_\kappa$ is nonincreasing; hence, $\ran f_\kappa$ is \kappaop.  Moreover, $\ran f_\kappa$ is dense in $P$.
\end{proof}

\begin{theorem}\label{THMcharpibounds}
For any space $X$ with point $p$, we have $\ochar{p,X}\leq\character{p,X}$ and $\opi{X}\leq\piweight{X}$ and $\ow{X}\leq\weight{X}^+$ and $\okchar{X}\leq\weight{X}$.
\end{theorem}
\begin{proof}
The first two relations immediately follow from Lemma~\ref{LEMkappasizeposet}; the third relation is trivial.  For the last relation, note that if $K$ is a compact subset of $X$, then it has neighborhood base of size at most $\weight{X}$; apply Lemma~\ref{LEMkappasizeposet}.
\end{proof}

Given Theorem~\ref{THMowprodsup}, justifying Observation~\ref{OBShomogbound} for $\ow{\cdot}$, $\opi{\cdot}$, and $\ochar{\cdot}$ amounts to justifying it for first countable homogeneous compacta, dyadic homogeneous compacta, and the two known kinds of exceptional homogeneous compacta.  The first countable case is the easiest.  By Arhangel${}^\prime$ski\u\i's Theorem, first countable compacta have weight at most $\cardc$, and therefore have Noetherian type at most $\cardc^+$.  Moreover, every point in a first countable space clearly has an \omegaop\ local base.  The only nontrivial bound is the one on Noetherian $\pi$\nbd-type.  For that, the following theorem suffices.

\begin{definition}\label{DEFpisw}
Give a space $X$, let $\pisw{X}$ denote the least $\kappa$ such that $X$ has a $\pi$\nbd-base $\mcA$ such that $\bigcap\mcB=\emptyset$ for all $\mcB\in[\mcA]^{\kappa^+}$.
\end{definition}

\begin{theorem}\label{THMopipisw}
If $X$ is a compactum, then $\opi{X}\leq\pisw{X}^+\leq\tightness{X}^+\leq\character{X}^+$.
\end{theorem}
\begin{proof}
Only the second relation is nontrivial; it is a theorem of Shapirovski\u\i~\cite{shapirovskii}.
\end{proof}

For dyadic homogeneous compacta, Theorem~\ref{THMdyadicmain} obviously implies Observation~\ref{OBShomogbound}; we will prove this theorem in Section~\ref{SECdyadic}.  Now consider the two known classes exceptional homogeneous compacta.  They are constructed by two techniques, resolutions and amalgams.  First we consider the exceptional resolution.

\begin{definition}\label{DEFresolution}
Suppose $X$ is a space, $\la Y_p\ra_{p\in X}$ is a sequence of nonempty spaces, and $\la f_p\ra_{p\in X}\in\prod_{p\in X}C(X\setminus\{p\},Y_p)$.  Then the \emph{resolution} $Z$ of $X$ at each point $p$ into $Y_p$ by $f_p$ is defined by setting $Z=\bigcup_{p\in X}\{p\}\times Y_p$ and declaring $Z$ to have weakest topology such that, for every $p\in X$, open neighborhood $U$ of $p$ in $X$, and open $V\subseteq Y_p$, the set $U\otimes V$ is open in $Z$ where $U\otimes V=(\{p\}\times V)\cup\bigcup_{q\in U\cap\inv{f}_p V}\{q\}\times Y_q$.
\end{definition}

The resolution of concern to us in constructed by van Mill~\cite{vanmill}.  It is a compactum with weight $\cardc$, $\pi$\nbd-weight $\omega$, and character $\omega_1$.  Moreover, assuming $\text{MA}+\neg\text{CH}$ (or just $\cardp>\omega_1$), this space is homogeneous.  (It is not homogeneous if $2^\omega<2^{\omega_1}$.)  For a proof that this space is exceptional (assuming $\text{MA}+\neg\text{CH}$), see \cite{milovicha}.  Clearly, this space has sufficiently small Noetherian type and $\pi$\nbd-type.  We just need to show that it has local Noetherian type $\omega$.  Van Mill's space is a resolution of $2^\omega$ at each point into $\mbT^{\omega_1}$ where $\mbT$ is the circle group $\reals/\integers$.

Notice that $\mbT$ is metrizable.  The following lemma proves that every metric compactum has Noetherian type $\omega$, along with some results that will be useful in Section~\ref{SECdyadic}.

\begin{lemma}\label{LEMcpctmetricbase}
Let $X$ be a metric compactum with base $\mcA$.  Then there exists $\mcB\subseteq\mcA$ satisfying the following.
\begin{enumerate}
\item\label{enumMetricbase} $\mcB$ is a base of $X$.
\item\label{enumMetricomegaop} $\mcB$ is \omegaop.
\item\label{enumMetricclose} If $U,V\in\mcB$ and $U\subsetneq V$, then $\closure{U}\subseteq V$.
\item\label{enumMetricfin} For all $\Gamma\in[\mcB]^{<\omega}$, there are only finitely many $U\in\mcB$ such that $\Gamma$ contains $\{V\in\mcB:U\subsetneq V\}$.
\end{enumerate}
\end{lemma}
\begin{proof}
Construct a sequence $\la\mcB_n\ra_{n<\omega}$ of finite subsets of $\mcA$ as follows.  For each $n<\omega$, let $E_n$ be the union of the set of all singletons in $\bigcup_{m<n}\mcB_m$.  Let $\mcC_n$ be the set of all $U\in\mcA$ for which and $U\cap E_n=\emptyset$ and
\begin{equation*}
2^{-n}\geq\diam U<\min\Biggl\{\diam V:V\in\bigcup_{m<n}\mcB_m\text{ and }0<\diam V\Biggr\}
\end{equation*}
and $\closure{U}\subseteq V$ for all $V\in\bigcup_{m<n}\mcB_m$ strictly containing $U$.  Then $\bigcup\mcC_n=X\setminus E_n$.  Let $\mcB_n$ be a minimal finite subcover of $\mcC_n$.  Set $\mcB=\bigcup_{n<\omega}\mcB_n$. To prove (\ref{enumMetricclose}, suppose $U\in\mcB_n$ and $V\in\mcB_m$ and $U\subsetneq V$.  Then $m\not=n$ by minimality of $\mcB_n$.  Also, $0<\diam V$ because $\emptyset\not=U\subsetneq V$.  Hence, if $m>n$, then $\diam V<\diam U$, in contradiction with $U\subsetneq V$.  Hence, $m<n$; hence, $\closure{U}\subseteq V$.

For (\ref{enumMetricbase}), let $p\in X$ and $n<\omega$, and let $V$ be the open ball with radius $2^{-n}$ and center $p$.  Then we just need to show that there exists $U\in\mcB$ such that $p\in U\subseteq V$.  Hence, \wma\ $\{p\}\not\in\mcB$.  Hence, $p\not\in E_{n+1}$; hence, there exists $U\in\mcB_{n+1}$ such that $p\in U$.  Since $\diam U\leq 2^{-n-1}$, we have $U\subseteq V$.  

For (\ref{enumMetricomegaop}), let $n<\omega$ and $U\in\mcB_n$.  If $U$ is a singleton, then every superset of $U$ in $\mcB$ is in $\bigcup_{m\leq n}\mcB_m$.  If $U$ is not a singleton, then $U$ has diamater at least $2^{-m}$ for some $m<\omega$; whence, every superset of $U$ in $\mcB$ is in $\bigcup_{l\leq m}\mcB_l$.

For (\ref{enumMetricfin}), suppose $\Gamma\in[\mcB]^{<\omega}$ and there exist infinitely many $U\in\mcB$ such that $\{V\in\mcB:U\subsetneq V\}\subseteq\Gamma$.  \Wma\ $\Gamma$ contains no singletons.  Choose an increasing sequence $\la k_n\ra_{n<\omega}$ in $\omega$ such that, for all $n<\omega$, there exists $U_n\in\mcB_{k_n}$ such that $\{V\in\mcB:U_n\subsetneq V\}\subseteq\Gamma$.  For each $n<\omega$, choose $p_n\in U_n$.  Since $\{U_n:n<\omega\}$ is infinite, we may choose $\la p_n\ra_{n<\omega}$ such that $\{p_n:n<\omega\}$ is infinite.  Let $p$ be an accumulation point of $\{p_n:n<\omega\}$.  Choose $m<\omega$ such that $2^{-m}<\diam V$ for all $V\in\Gamma$.  Since $p$ is not an isolated point, there exists $W\in\mcB_m$ such that $p\in W$.  Then $W\not\in\Gamma$; hence, $W$ does not strictly contain $U_n$ for any $n<\omega$.  Choose $q\in W\setminus\{p\}$ such that $W$ contains $\{x:d(p,x)\leq d(p,q)\}$; set $r=d(p,q)$.  Let $B$ be the open ball of radius $r/2$ centered about $p$.  Then there exists $n<\omega$ such that $2^{-k_n}<r/2$ and $p_n\in B$.  Hence, $\diam U_n<r/2$ and $U_n\cap B\not=\emptyset$; hence, $U_n\subseteq W$ and $q\not\in U_n$; hence, $U_n\subsetneq W$, which is absurd.  Therefore, for each $\Gamma\in[\mcB]^{<\omega}$, there are only finitely many $U\in\mcB$ such that $\{V\in\mcB:U\subsetneq V\}\subseteq\Gamma$.
\end{proof}

We have $\ow{2^\omega}=\ow{\mbT^{\omega_1}}=\omega$ by Lemma~\ref{LEMcpctmetricbase} and Theorem~\ref{THMowprodsup}.  Therefore, the following theorem implies that van Mill's space has local Noetherian type $\omega$.

\begin{lemma}[\cite{vanmill}]\label{LEMresollocbase}
Suppose $X$, $\la Y_p\ra_{p\in X}$, $\la f_p\ra_{p\in X}$, and $Z$ are as in Definition~\ref{DEFresolution}.  Suppose $\mcU$ is a local base at a point $p$ in $X$ and $\mcV$ is a local base at a point $y$ in $Y_p$.  Then $\{U\otimes V:\la U,V\ra\in\mcU\times(\mcV\cup\{Y_p\})\}$ is a local base at $\la p,y\ra$ in $Z$.
\end{lemma}

\begin{theorem}\label{THMocharres}
Suppose $X$, $\la Y_p\ra_{p\in X}$, $\la f_p\ra_{p\in X}$, and $Z$ are as in Definition~\ref{DEFresolution}.  Then $\ochar{\la p,y\ra,Z}\leq\ow{X}\ochar{y,Y_p}$ for all $\la p,y\ra\in Z$.
\end{theorem}
\begin{proof}
Set $\kappa=\ow{X}\ochar{y,Y_p}$.  Let $\mcA$ be a \kappaop\ base of $X$ and let $\mcB$ be a \kappaop\ local base at $y$ in $Y_p$; \wma\ $Y_p\in\mcB$.  Set $\mcC=\{U\in\mcA:p\in U\}$.  Set $\mcD=\{U\otimes V:\la U,V\ra\in\mcC\times\mcB\}$, which is a local base at $\la p,y\ra$ in $Z$ by Lemma~\ref{LEMresollocbase}.  If there exists $U\otimes V\in\mcD$ such that $U\cap\inv{f}_p V=\emptyset$, then $U\otimes V$ is homeomorphic to $V$; whence, $\ochar{\la p,y\ra,Z}=\ochar{y,Y_p}\leq\kappa$.  Hence, \wma\ $U\cap\inv{f}_p V\not=\emptyset$ for all $U\otimes V\in\mcD$.

It suffices to show that $\mcD$ is \kappaop.  Suppose $U_i\otimes V_i\in\mcD$ for all $i<2$ and $U_0\otimes V_0\subseteq U_1\otimes V_1$.  Then $V_0\subseteq V_1$ and $\emptyset\not=U_0\cap\inv{f}_p V_0\subseteq U_1\cap\inv{f}_p V_1$.  Since $\mcB$ is \kappaop, there are fewer than $\kappa$\nbd-many possibilities for $V_1$ given $V_0$.  Since $\mcA$ is a \kappaop\ base, there are fewer than $\kappa$\nbd-many possibilities for $U_1$ given $U_0$ and $V_0$.  Hence, there are fewer than $\kappa$\nbd-many possibilities for $U_1\otimes V_1$ given $U_0\otimes V_0$.
\end{proof}

\begin{definition}
Let $\cardp$ denote the least $\kappa$ for which some $\mcA\in[[\omega]^\omega]^\kappa$ has the strong finite intersection property but does not have a nontrivial pseudointersection.  By a theorem of Bell\cite{bell}, $\cardp$ is also the least $\kappa$ for which there exist a $\sigma$\nbd-centered poset $\mbP$ and a family $\mcD$ of $\kappa$\nbd-many dense subsets of $\mbP$ such that $\mbP$ does not have a $\mcD$\nbd-generic filter.
\end{definition}

\begin{definition}
Given a space $X$, let $\Aut{X}$ denote the set of its autohomeomorphisms.
\end{definition}

Van Mill's construction has been generalized by Hart and Ridderbos\cite{hart}.  They show that one can produce an exceptional homogeneous compactum with weight $\cardc$ and $\pi$\nbd-weight $\omega$ by carefully resolving each point of $2^\omega$ into a fixed space $Y$ satisfying the following conditions.
\begin{enumerate}
\item\label{enumYcpct} $Y$ is a homogeneous compactum.
\item\label{enumYdwp} $\omega_1\leq\character{Y}\leq\weight{Y}<\cardp$.
\item\label{enumYdenseorbit} $\exists d\in Y\ \ \exists\eta\in\Aut{Y}\ \ \closure{\{\eta^n(d):n<\omega\}}=Y$.
\item If $\gamma\omega$ is a compactification of $\omega$ and $\gamma\omega\setminus\omega\homeo Y$, then $Y$ is a retract of $\gamma\omega$.
\end{enumerate}
By Theorem~\ref{THMocharres}, to show that such resolutions have local Noetherian type $\omega$, it suffices to show that every such $Y$ has local Noetherian type $\omega$.  Theorem~\ref{THMhomogwpochar} will accomplish this.

\begin{theorem}\label{THMocharpicharchar}
Suppose $X$ is a compactum and $\pi\character{p,X}=\character{q,X}$ for all $p,q\in X$.  Then $\ochar{p,X}=\omega$ for some $p\in X$.  In particular, if $X$ is a homogeneous compactum and $\pi\character{X}=\character{X}$, then $\ochar{X}=\omega$.
\end{theorem}

The proof of Theorem~\ref{THMocharpicharchar} will be delayed until Section~\ref{SEChomogloc}.

The following lemma is essentially a generalization of a similar result of Juh\'asz\cite{juhasz89}.

\begin{lemma}\label{LEMsepcpctwcardp}
Suppose $X$ is a compactum and $\omega=\density{X}\leq\weight{X}<\cardp$.  Then there exists $p\in X$ such that $\character{p,X}\leq\piweight{X}$.
\end{lemma}
\begin{proof}
Let $\mcA$ be a base of $X$ of size at most $\weight{X}$.  Let $\mcB$ be a $\pi$\nbd-base of $X$ of size at most $\piweight{X}$.  For each $\la U,V\ra\in\mcB^2$ satisfying $\closure{U}\subseteq V$, choose a closed $G_\delta$\nbd-set $\Phi(U,V)$ such that $\closure{U}\subseteq\Phi(U,V)\subseteq V$.  Then $\ran\Phi$, ordered by $\subseteq$, is $\sigma$\nbd-centered because $\density{X}=\omega$.  Since $\card{\mcA}<\cardp$, there is a filter $\mcG$ of $\ran\Phi$ such that for all disjoint $U,V\in\mcA$ some $K\in\mcG$ satisfies $U\cap K=\emptyset$ or $V\cap K=\emptyset$.  Hence, there exists a unique $p\in\bigcap\mcG$.  Hence, $p$ has pseudocharacter, and therefore character, at most $\card{\mcG}$, which is at most $\piweight{X}$.
\end{proof}

\begin{theorem}\label{THMhomogwpochar}
If $X$ is a homogeneous compactum and $\omega=\density{X}\leq\weight{X}<\cardp$, then $\ochar{X}=\omega$.
\end{theorem}
\begin{proof}
By Lemma~\ref{LEMsepcpctwcardp}, $\character{X}\leq\piweight{X}=\pi\character{X}\density{X}=\pi\character{X}$.  Hence, by Theorem~\ref{THMocharpicharchar}, $\ochar{X}=\omega$.
\end{proof}

Amalgams are defined in \cite{milovicha} as follows.

\begin{definition}\label{DEFamalgam}
Suppose $X$ is a $T_0$ space, $\msS$ is a subbase of $X$ such that $\emptyset\not\in\msS$, and $\la Y_S\ra_{S\in\msS}$ is a sequence of nonempty spaces.  
The \emph{amalgam} $Y$ of $\la Y_S: {S\in\msS}\ra$ is defined by setting
$Y=\bigcup_{p\in X}\prod_{p\in S\in\msS} Y_S$
and declaring $Y$ to have the weakest topology such that, for each $S\in\msS$ and open $U\subseteq Y_S$, the set $\inv{\pi}_S U$ is open in $Y$ where $\inv{\pi}_S U=\{p\in Y: S\in\dom p\text{ and }p(S)\in U\}$.  Define $\pi:Y\rightarrow X$ by $\{\pi(p)\}=\bigcap\dom p$ for all $p\in Y$.  It is easily verified that $\pi$ is continuous.
\end{definition}  

\begin{theorem}\label{THMamalgam}
Suppose $X$, $\msS$, $\la Y_S\ra_{S\in\msS}$, and $Y$ be as in Definition~\ref{DEFamalgam}.  Then we have the following relations for all $p\in Y$.
\begin{align*}
\ow{Y} &\leq\ow{X}\sup_{S\in\msS}\ow{Y_S}\\
\opi{Y} &\leq\opi{X}\sup_{S\in\msS}\opi{Y_S}\\
\ochar{p,Y} &\leq\ochar{\pi(p),X}\sup_{S\in\dom p}\ochar{p(S),Y_S}
\end{align*}
\end{theorem}
\begin{proof}
We will only prove the first relation; the proofs of the others are almost identical.  Set $\kappa=\ow{X}\sup_{S\in\msS}\ow{Y_S}$.  Let $\mcA$ be a \kappaop\ base of $X$.  For each $S\in\msS$, let $\mcB_S$ be a \kappaop\ base of $Y_S$.    Set 
\begin{equation*}
\mcC=\biggl\{\inv{\pi} U\cap\bigcap_{S\in\dom\tau}\inv{\pi}_S \tau(S):\tau\in\bigcup_{\mcF\in[\msS]^{<\omega}}\prod_{S\in\mcF}\mcB_S\setminus\{Y_S\}\text{ and } \mcA\ni U\subseteq\bigcap\dom\tau\biggr\}.
\end{equation*}
Then $\mcC$ is clearly a base of $Y$.  Let us show that $\mcC$ is \kappaop.  Suppose $\inv{\pi} U_i\cap\bigcap_{S\in\dom\tau_i}\inv{\pi}_S\tau_i(S)\in\mcC$ for all $i<2$ and
\begin{equation*}
\inv{\pi} U_0\cap\bigcap_{S\in\dom\tau_0}\inv{\pi}_S\tau_0(S)\subseteq\inv{\pi} U_1\cap\bigcap_{S\in\dom\tau_1}\inv{\pi}_S\tau_1(S).
\end{equation*}
Then $U_0\subseteq U_1$ and $\dom\tau_0\supseteq\dom\tau_1$ and $\tau_0(S)\subseteq\tau_1(S)$ for all $S\in\dom\tau_1$.  Hence, there are fewer than $\kappa$\nbd-many possibilities for $U_1$ and $\tau_1$ given $U_0$ and $\tau_0$.
\end{proof}

An exceptional homogeneous compactum $Y$ is constructed in \cite{milovicha} with $X=\mbT$ and $\weight{Y_S}=\piweight{Y_S}=\cardc$ and $\character{Y_S}=\omega$ for all $S\in\msS$.  Hence, $\ow{Y_S}\leq\cardc^+$ and $\ochar{Y_S}=\omega$ for each $S\in\msS$.  Moreover, each $Y_S$ is $2^\gamma$ ordered lexicographically where $\gamma$ is a fixed indecomposable ordinal in $\omega_1\setminus(\omega+1)$.  Since $\cf\gamma=\omega$, it is easy to construct an \omegaop\ $\pi$\nbd-base of this space.  Hence, by Theorem~\ref{THMamalgam}, $\ow{Y}\leq\cardc^+$ and $\opi{Y}=\ochar{Y}=\omega$.  Thus, Observation~\ref{OBShomogbound} is justified for $\ow{\cdot}$, $\opi{\cdot}$, and $\ochar{\cdot}$.

It remains to justify Observation~\ref{OBShomogbound} for $\okchar{\cdot}$.  We first note that all known homogeneous compacta are continuous images of products of compacta each of weight at most $\cardc$.  (Moreover, it it shown in \cite{milovicha} that any $Z$ as in Definition~\ref{DEFamalgam} is a continuous image of $X\times\prod_{S\in\msS}Y_S$.)  Therefore, the following theorem will suffice.

\begin{theorem}\label{THMokcharctsprodw}
Suppose $Y$ is a continuous image of a product $X=\prod_{i\in I}X_i$ of compacta.  Then $\okchar{Y}\leq\sup_{i\in I}\weight{X_i}$
\end{theorem}

Before proving the above theorem, we first prove two lemmas.

\begin{definition}
Given subsets $P$ and $Q$ of a common poset, define $P$ and $Q$ to be \emph{mutually dense} if for all $p_0\in P$ and $q_0\in Q$ there exist $p_1\in P$ and $q_1\in Q$ such that $p_0\geq q_1$ and $q_0\geq p_1$.
\end{definition}

\begin{lemma}\label{LEMmutuallydense}
Let $\kappa$ be a cardinal and let $P$ and $Q$ be mutually dense subsets of a common poset.  Then $P$ is almost \kappaop\ if and only if $Q$ is.
\end{lemma}
\begin{proof}
Suppose $D$ is a \kappaop\ dense subset of $P$.  Then it suffices to construct a \kappaop\ dense subset of $Q$.  Define a partial map $f$ from $\card{D}^+$ to $Q$ as follows.  Set $f_0=\emptyset$.  Suppose $\alpha<\card{D}^+$ and we have constructed a partial map $f_\alpha$ from $\alpha$ to $Q$.  Set $E=\{d\in D: d\not\geq q\text{ for all }q\in\ran f_\alpha\}$.  If $E=\emptyset$, then set $f_{\alpha+1}=f_\alpha$.  Otherwise, choose $q\in Q$ such that $q\leq e$ for some $e\in E$, and let $f_{\alpha+1}$ be the smallest function extending $f_\alpha$ such that $f_{\alpha+1}(\alpha)=q$.  For limit ordinals $\gamma\leq\card{D}^+$, set $f_\gamma=\bigcup_{\alpha<\gamma}f_\alpha$.  Set $f=f_{\card{D}^+}$.

Let us show that $\ran f$ is \kappaop.  Suppose otherwise.  Then there exists $q\in\ran f$ and an increasing sequence $\la\xi_\alpha\ra_{\alpha<\kappa}$ in $\dom f$ such that $q\leq f(\xi_\alpha)$ for all $\alpha<\kappa$.  By the way we constructed $f$, there exists $\la d_\alpha\ra_{\alpha<\kappa}\in D^\kappa$ such that $f(\xi_\beta)\leq d_\beta\not=d_\alpha$ for all $\alpha<\beta<\kappa$.  Choose $p\in P$ such that $p\leq q$.  Then choose $d\in D$ such that $d\leq p$.  Then $d\leq d_\beta\not=d_\alpha$ for all $\alpha<\beta<\kappa$, which contradicts that $D$ is \kappaop.  Therefore, $\ran f$ is \kappaop.

Finally, let us show that $\ran f$ is a dense subset of $Q$.  Suppose $q\in Q$.  Choose $p\in P$ such that $p\leq q$.  Then choose $d\in D$ such that $d\leq p$.  By the way we constructed $f$, there exists $r\in\ran f$ such that $r\leq d$; hence, $r\leq q$.
\end{proof}

\begin{lemma}\label{LEMokcharcts}
Suppose $f\colon X\rightarrow Y$ is a continuous surjection between compacta and $C$ is closed in $Y$.  Then $\ochar{\inv{f}C,X}=\ochar{C,Y}$.
\end{lemma}
\begin{proof}
Let $\mcA$ be a neighborhood base of $C$.  By Lemma~\ref{LEMmutuallydense}, it suffices to show that $\{\inv{f}V:V\in\mcA\}$ is a neighborhood base of $\inv{f}C$.  Suppose $U$ is a neighborhood of $\inv{f}C$.  By normality of $Y$, we have $\inv{f}C=\bigcap_{V\in\mcA}\inv{f}\closure{V}$.  By compactness of $X$, we have $\inv{f}\closure{V}\subseteq U$ for some $V\in\mcA$.  Thus, $\{\inv{f}V:V\in\mcA\}$ is a neighborhood base of $\inv{f}C$ as desired. 
\end{proof}

\begin{proof}[Proof of Theorem~\ref{THMokcharctsprodw}]
By Lemma~\ref{LEMokcharcts}, \wma\ $Y=X$.  By Theorem~\ref{THMowprodsup}, \wma\ $I$ is finite.  Apply Theorem~\ref{THMcharpibounds}.
\end{proof}

How sharp are the bounds of Observation~\ref{OBShomogbound}?  (\ref{enumobsochar}) is trivially sharp as every space has local Noetherian type at least $\omega$.  We will show that there is a homogeneous compactum with Noethian type $\cardc^+$, namely, the double arrow space.  Moreover, we will show that Suslin lines have uncountable Noetherian $\pi$\nbd-type.   It is known to be consistent that there are homogeneous compact Suslin lines, but it is also known to be consistent that there are no Suslin lines.
It is not clear whether it is consistent that all homogeneous compacta have Noetherian $\pi$\nbd-type $\omega$, even if we restrict to the first countable case.  Also, it is not clear in any model of ZFC whether all first countable homogeneous compacta have compact Noetherian type $\omega$.

\begin{question}
Is there a first countable compactum with uncountable compact Noetherian type?
\end{question}

The following proposition is essentially due to Peregudov~\cite{peregudov97}.

\begin{proposition}\label{PROowpiw}
If $X$ is a space and $\piweight{X}<\cf\kappa\leq\kappa\leq\weight{X}$, then $\ow{X}>\kappa$.
\end{proposition}
\begin{proof}
Suppose $\mcA$ is a base of $X$ and $\mcB$ is $\pi$\nbd-base of $X$ of size $\piweight{X}$.  Then $\card{\mcA}\geq\kappa$; hence, there exist $\mcU\in[\mcA]^\kappa$ and $V\in\mcB$ such that $V\subseteq\bigcap\mcU$.  Hence, there exists $W\in\mcA$ such that $W\subseteq V\subseteq\bigcap\mcU$; hence, $\mcA$ is not \kappaop.
\end{proof}

\begin{example}
The double arrow space, defined as $((0,1]\times\{0\})\cup([0,1)\times\{1\})$ ordered lexicographically, has $\pi$\nbd-weight $\omega$ and weight $\cardc$, and is is known to be compact and homogeneous.  By Proposition~\ref{PROowpiw}, it has Noetherian type $\cardc^+$.
\end{example}

\begin{theorem}
Suppose $X$ is a Suslin line.  Then $\opi{X}\geq\omega_1$.
\end{theorem}
\begin{proof}
Let $\mcA$ be a $\pi$\nbd-base of $X$ consisting only of open intervals.  By Lemma~\ref{LEMmutuallydense}, it suffices to show that $\mcA$ is not \omegaop.  Construct a sequence $\la\mcB_n\ra_{n<\omega}$ of maximal pairwise disjoint subsets of $\mcA$ as follows.  Choose $\mcB_0$ arbitrarily.  Given $n<\omega$ and $\mcB_n$, choose $\mcB_{n+1}$ such that it refines $\mcB_n$ and $\mcB_n\cap\mcB_{n+1}\subseteq[X]^1$.

Let $E$ denote the set of all endpoints of intervals in $\bigcup_{n<\omega}\mcB_n$.  Since $X$ is Suslin, there exists $U\in\mcA\setminus[X]^1$ such that $U\cap E=\emptyset$.  For each $n<\omega$, the set $\bigcup\mcB_n$ is dense in $X$ by maximality; whence, there exists $V_n\in\mcB_n$ such that $U\cap V_n\not=\emptyset$.  Since $U\cap E=\emptyset$, we have $U\subseteq\bigcap_{n<\omega}V_n$.  Thus, $\mcA$ is not \omegaop.
\end{proof}

$\mathrm{MA}+\neg\mathrm{CH}$ implies there are no Souslin lines.  It is not clear whether it further implies every homogeneous compactum has Noetherian $\pi$\nbd-type $\omega$.  However, the next theorem gives us a partial result.  First, we need a lemma very similar to the result that $\mathrm{MA}+\neg\mathrm{CH}$ implies all Aronszajn trees are special.

\begin{definition}
Given a subset $E$ of a poset $Q$, let $\upset_Q E$ denote the set of $q\in Q$ for which $q$ has a lower bound in $E$.
\end{definition}

\begin{lemma}\label{LEMgenspecialaronszajn}
Assume MA.  Suppose $Q$ is an \cardop{\omega_1}\ poset of size less than $\cardc$.  Then $Q$ is almost \omegaop\ or $Q$ has an uncountable centered subset.
\end{lemma}
\begin{proof}
Set $\mbP=[Q]^{<\omega}$ and order $\mbP$ such that $\sigma\leq\tau$ if and only if $\sigma\cap\upset_Q\tau=\tau$.  A sufficiently generic filter $G$ of $\mbP$ will be such that $\bigcup G$ is a dense \omegaop\ subset of $Q$.  Hence, if $\mbP$ is ccc, then $Q$ is almost \omegaop.  Hence, \wma\ $\mbP$ has an antichain $A$ of size $\omega_1$.  \Wma\ $A$ is a $\Delta$\nbd-system with root $\rho$.  Since $Q$ is \cardop{\omega_1}, \wma\ $\sigma\cap\upset_Q\rho=\rho$ for all $\sigma\in A$.  Choose a bijection $\la a_\alpha\ra_{\alpha<\omega_1}$ from $\omega_1$ to $A$.  \Wma\ there exists an $n<\omega$ such that $\card{a_\alpha\setminus\rho}=n$ for all $\alpha<\omega_1$.  For each $\alpha<\omega_1$, choose a bijection $\la a_{\alpha,i}\ra_{i<n}$ from $n$ to $a_\alpha\setminus\rho$.  For each $x\in Q$ and $i<n$, set $E_{x,i}=\{\alpha<\omega_1:x\leq_Q a_{\alpha,i}\text{ or }a_{\alpha,i}\leq_Q x\}$.  For each $\alpha<\omega_1$, since $A$ is an antichain, we have $\bigcup_{i<n}\bigcup_{j<n}E_{a_{\alpha,i},j}=\omega_1$.  Choose a uniform ultrafilter $\mcU$ on $\omega_1$.  Then we may choose $B\in[(\bigcup A)\setminus\rho]^{\omega_1}$ and $i<n$ such that $E_{x,i}\in\mcU$ for all $x\in B$.  

It suffices to show that $B$ is centered.  Let $\sigma\in[B]^{<\omega}$.  Set $E=\bigcap_{x\in\sigma}E_{x,i}$.  Then $E\in\mcU$; hence, $\card{E}=\omega_1$; hence, we may choose $\alpha\in E\setminus\{\beta<\omega_1:a_{\beta,i}\in\upset_Q\sigma\}$.  Then $a_{\alpha,i}<_Q x$ for all $x\in\sigma$.  Thus, $B$ is centered.
\end{proof}

\begin{lemma}\label{LEMirredopi}
Suppose $f\colon X\rightarrow Y$ is an irreducible continuous surjection between spaces.  Then $\opi{X}=\opi{Y}$.
\end{lemma}
\begin{proof}
Let $\mcA$ be a \cardop{\opi{X}}\ $\pi$\nbd-base of $X$ and let $\mcB$ be a \cardop{\opi{Y}}\ $\pi$\nbd-base of $Y$.  By Lemma~\ref{LEMmutuallydense}, \wma\ $\mcA$ consists only of regular open sets.  Set $\mcC=\{\inv{f}U:U\in\mcB\}$.  Then $\mcC$ is \cardop{\opi{Y}}.  Suppose $U$ is a nonempty open subset of $X$.  Then we may choose $V\in\mcB$ such that $V\cap f[X\setminus U]=\emptyset$.  Then $\inv{f}V\subseteq U$.  Thus, $\mcC$ is a $\pi$\nbd-base of $X$; hence, $\opi{X}\leq\opi{Y}$.

Set $\mcD=\{Y\setminus f[X\setminus U]:U\in\mcA\}$.  Suppose $V$ is a nonempty open subset of $Y$.  Then we may choose $U\in\mcA$ such that $U\subseteq\inv{f}V$.  Then $Y\setminus f[X\setminus U]\subseteq V$.  Thus, $\mcD$ is a $\pi$\nbd-base of $Y$.  Now suppose $U_0,U_1\in\mcA$ and $U_0\not\subseteq U_1$.  Then $U_0\not\subseteq \closure{U}_1$ by regularity.  By irreducibility, we may choose $p\in Y\setminus f[X\setminus(U_0\setminus\closure{U}_1)]$.  Then $p\in f[X\setminus U_1]$ and $p\not\in f[X\setminus U_0]$.  Hence, $Y\setminus f[X\setminus U_0]\not\subseteq Y\setminus f[X\setminus U_1]$.  Thus, $\mcD$ is \cardop{\opi{X}}; hence, $\opi{Y}\leq\opi{X}$.
\end{proof}

\begin{theorem}
Assume MA.  Let $X$ be a compactum such that $\tightness{X}=\omega$ and $\piweight{X}<\cardc$.  Then $\opi{X}=\omega$.
\end{theorem}
\begin{proof}
\Wma\ $X$ is a closed subspace of $[0,1]^\kappa$ for some cardinal $\kappa$.
By a result of Shapirovski\u\i~\cite{shapirovskii},
since $\tightness{X}=\omega$, there is an irreducible continuous map $f$ from $X$ onto a subspace of $\bigcup_{I\in[\kappa]^\omega}[0,1]^I\times\{0\}^{\kappa\setminus I}$.  By Lemma~\ref{LEMirredopi}, \wma\ $X\subseteq\bigcup_{I\in[\kappa]^\omega}[0,1]^I\times\{0\}^{\kappa\setminus I}$.  Set $\mcF=\Fn{\kappa}{(\mbQ\cap(0,1])^2}$ and 
\begin{equation*}
 \mcA=\left\{X\cap\bigcap_{\alpha\in\dom\sigma}\inv{\pi_\alpha}\left(\sigma(\alpha)(0),\sigma(\alpha)(1)\right):\sigma\in\mcF\right\}\setminus\{\emptyset\},
\end{equation*}
which is a $\pi$\nbd-base of $X$.  By Theorem~\ref{THMopipisw} and Lemma~\ref{LEMmutuallydense}, $\mcA$ contains an \cardop{\omega_1}\ dense subset $\mcB$, and it suffices to show that $\mcB$ is almost \omegaop.  Seeking a contradiction, suppose $\mcB$ is not almost \omegaop.  By Lemma~\ref{LEMgenspecialaronszajn}, $\mcB$ contains an uncountable centered subset $\mcC$.  Let the map
\begin{equation*}
\left\la X\cap\bigcap_{\alpha\in\dom\sigma_\beta}\inv{\pi_\alpha}(\sigma_\beta(\alpha)(0),\sigma_\beta(\alpha)(1))\right\ra_{\beta<\omega_1}
\end{equation*}
be an injection from $\omega_1$ to $\mcC$.  Then $\card{\bigcup_{\beta<\omega_1}\dom\sigma_\beta}=\omega_1$.  By compactness, the set
\begin{equation*}
X\cap\bigcap_{\beta<\omega_1}\bigcap_{\alpha\in\dom\sigma_\beta}\inv{\pi_\alpha}[\sigma_\beta(\alpha)(0),\sigma_\beta(\alpha)(1)]
\end{equation*}
is nonempty, in contradiction with $X\subseteq\bigcup_{I\in[\kappa]^\omega}[0,1]^I\times\{0\}^{\kappa\setminus I}$.  
\end{proof}

Concerning compact Noetherian type, we note that if there is a homogeneous compactum $X$ for which $\okchar{X}\geq\omega_1$, then $X$ is not an ordered space.

\begin{definition}
A point $p$ in a space $X$ is $P_\kappa$\nbd-\emph{point} if, for every set $\mcA$ of fewer than $\kappa$\nbd-many neighborhoods of $p$, the set $\bigcap\mcA$ has $p$ in its interior.  A $P$\nbd-point is a $P_{\omega_1}$\nbd-point.
\end{definition}

\begin{theorem}\label{THMhomogordokchar}
If $X$ is a homogeneous ordered compactum, then $\okchar{X}=\omega$.
\end{theorem}
\begin{proof}
\Wma\ $X$ is infinite; hence, $X$ has a point that is not a $P$\nbd-point.  By homogeneity, $\min X$ is not a $P$\nbd-point; hence, $\min X$ has countable character.  By homogeneity, $X$ is first countable.  Let $C$ be closed in $X$.  Then $X\setminus C$ is a disjoint union of open intervals $\bigcup_{i\in I}(a_i,b_i)$ such that $(a_i,b_i)=\bigcup_{n<\omega}[a_{i,n},b_{i,n}]$ and $\la a_{i,n}\ra_{n<\omega}$ is nonincreasing and $\la b_{i,n}\ra_{n<\omega}$ is nondecreasing for all $i\in I$.  Hence, $\{X\setminus\bigcup_{i\in\dom\sigma}[a_{i,\sigma(i)},b_{i,\sigma(i)}]:\sigma\in\Fn{I}{\omega}\}$ is an \omegaop\ neighborhood base of $C$.
\end{proof}

It is worth noting that while products do not decrease cellularity, they can decrease $\ow{\cdot}$, $\opi{\cdot}$, and $\ochar{\cdot}$, as shown by the following theorem of Malykhin~\cite{malykhin}.

\begin{theorem}\label{THMowproductweight}
Let $p\in X=\prod_{i\in I}X_i$ where $X_i$ is a nonsingleton $T_1$ space for all $i\in I$.  If $\sup_{i\in I}\weight{X_i}\leq\card{I}$, then $\ow{X}=\omega$.  If $\sup_{i\in I}\piweight{X_i}\leq\card{I}$, then $\opi{X}=\omega$.  If $\sup_{i\in I}\character{p(i),X_i}\leq\card{I}$, then $\ochar{p,X}=\omega$.
\end{theorem}
\begin{proof}
See \cite{malykhin} for a proof of the first implication.  That proof can be easily modified to demonstrate the other implications.
\end{proof}

In constrast, $\okchar{\cdot}$ is not decreased by products when the factors are compacta.  Just as is true of cellularity, the compact Noetherian type of a product of compacta is the supremum of the compact Noetherian types of its finite subproducts.

\begin{theorem}\label{THMokcharprod}
If $X=\prod_{i\in I}X_i$ is a product of compacta, then $\okchar{X}=\sup_{\sigma\in[I]^{<\omega}}\okchar{\prod_{i\in\sigma}X_i}$.
\end{theorem}
\begin{proof}
To prove ``$\leq$'', apply Theorem~\ref{THMowprodsup}.  To prove ``$\geq$'', apply Lemma~\ref{LEMokcharcts}.
\end{proof}

Though cellularity and compact Noetherian type behave similarly for compacta, they do not coincide, even assuming homogeneity.  Given any indecomposable ordinal $\gamma$ strictly between $\omega$ and $\omega_1$, the lexicographic ordering of $2^\gamma$ is homogeneous and compact and has cellularity $\cardc$ by a result of Maurice~\cite{maurice}.  However, by Theorem~\ref{THMhomogordokchar}, this space has compact Noetherian type $\omega$.

\section{Dyadic compacta}\label{SECdyadic}

In this section, we prove a strengthened version of Theorem~\ref{THMdyadicmain} and generalize it to continuous images of products of compacta with bounded weight.  We also investigate the spectrum of Noetherian types of dyadic compacta.  Our approach is to start with results about subsets of free boolean algebras and then use Stone duality to apply them to families of open subsets of dyadic compacta.

By Lemma~\ref{LEMkappasizeposet}, every countable subset of a free boolean algebra is almost \omegaop.  We wish to prove  this for all subsets of free boolean algebras.  We achieve this by approximating free boolean algebras by smaller free subalgebras using elementary substructures.  More specifically, we use elementary submodels of $H_\theta$ where $\theta$ is a regular cardinal and $H_\theta$ is the $\{\in\}$\nbd-structure of the family of sets that hereditarily have size less than $\theta$.  Whenever we use $H_\theta$ in an argument, we implicitly assume that $\theta$ is sufficiently large to make the argument valid.  As is typical with elementary submodels of $H_\theta$, we need reflection properties.  For our purposes, the crucial reflection property of free boolean algebras is given by the following lemma.

\begin{lemma}\label{LEMfreeboolreflect}
Let $B$ be a free boolean algebra and let $\{B,\land,\lor\}\subseteq M\elemsub H_\theta$.  Then, for all $q\in B$, there exists $r\in B\cap M$ such that, for all $p\in B\cap M$, we have $p\geq q$ if and only if $p\geq r$.  In particular, $r\geq q$.
\end{lemma}
\begin{proof}
Let $q\in B$.  \Wma\ $q\not=0$.  By elementarity, there exists a map $g\in M$ enumerating a set of mutually independent generators of $B$.  Set  $G=\bigcup\{\{g(i),g(i)'\}:i\in\dom g\}$.  Then there exists $\eta\in[[G]^{<\omega}]^{<\omega}$ such that $q=\bigvee_{\tau\in\eta}\bigwedge\tau$ and $\bigwedge\tau\not=0$ for all $\tau\in\eta$.  Set $r=\bigvee_{\tau\in\eta}\bigwedge(\tau\cap M)$.  Let $p\in B\cap M$; \wma\ $p\not=1$.  Then there exists $\zeta\in[[G\cap M]^{<\omega}]^{<\omega}$ such that $p=\bigwedge_{\sigma\in\zeta}\bigvee\sigma$ and $\bigvee\sigma\not=1$ for all $\sigma\in\zeta$.  Hence, $p\geq q$ iff, for all $\sigma\in\zeta$ and $\tau\in\eta$, we have $\bigvee\sigma\geq\bigwedge\tau$, which is equivalent to $\sigma\cap\tau\not=\emptyset$, which is equivalent to $\sigma\cap\tau\cap M\not=\emptyset$.  Thus, $p\geq q$ if and only if $p\geq r$.
\end{proof}

\begin{theorem}\label{THMomegaopbool}
Every subset of every free boolean algebra is almost \omegaop.
\end{theorem}
\begin{proof}
Let $B$ be a free boolean algebra; set $\kappa=\card{B}$.  Given $A\subseteq B$, let $\upset A$ denote the smallest semifilter of $B$ containing $A$; if $A=\{a\}$ for some $a$, then set $\upset a=\upset A$.  Let $Q$ be a subset of $B$.  If $Q$ is a countable, then $Q$ is almost \omegaop\ by Lemma~\ref{LEMkappasizeposet}.  Therefore, \wma\ that $\kappa>\omega$ and the theorem is true for all free boolean algebras of size less than $\kappa$.

We will construct a continuous elementary chain $\la M_\alpha\ra_{\alpha<\kappa}$ of elementary submodels of $H_\theta$ and a continuous increasing sequence of sets $\la D_\alpha\ra_{\alpha<\kappa}$ satisfying the following conditions for all $\alpha<\kappa$.
\begin{enumerate}
\item\label{enumMcontains} $\alpha\cup\{B,\land,\lor,Q\}\subseteq M_\alpha$ and $\card{M_\alpha}\leq\card{\alpha}+\omega$.
\item\label{enumDdense} $D_\alpha$ is a dense subset of $Q\cap M_\alpha$.
\item\label{enumDomegaop} $D_\alpha\cap\upset q$ is finite for all $q\in Q\cap M_\alpha$.
\item\label{enumDstable} $D_{\alpha+1}\cap\upset q=D_\alpha\cap\upset q$ for all $q\in Q\cap M_\alpha$.
\end{enumerate}
Given this construction, set $D=\bigcup_{\alpha<\kappa}D_\alpha$.  Then $D$ is a dense subset of $Q$ by (\ref{enumDdense}).  Moreover, if $\alpha<\kappa$ and $d\in D_\alpha$, then $d\in Q\cap M_\alpha$ by (\ref{enumDdense}); whence, $d$ is below at most finitely many elements of $D$ by (\ref{enumDomegaop}) and (\ref{enumDstable}).  Hence, $Q$ is almost \omegaop.

For stage $0$, choose any $M_0\elemsub H_\theta$ satisfying (\ref{enumMcontains}).  Since $Q\cap M_0\subseteq B\cap M_0$, we may choose $D_0$ to be an \omegaop\ dense subset of $Q\cap M_0$, exactly what (\ref{enumDdense}) and (\ref{enumDomegaop}) require.  At limit stages, (\ref{enumMcontains}) and (\ref{enumDdense}) are clearly preserved, and (\ref{enumDomegaop}) is preserved because of (\ref{enumDstable}).

For a successor stage $\alpha+1$, choose $M_{\alpha+1}$ such that $M_\alpha\elemsub M_{\alpha+1}\elemsub H_\theta$ and (\ref{enumMcontains}) holds for stage $\alpha+1$.  Since $Q\cap M_{\alpha+1}\subseteq B\cap M_{\alpha+1}$, there is an \omegaop\ dense subset $E$ of $Q\cap M_{\alpha+1}$.  Set $D_{\alpha+1}=D_\alpha\cup(E\setminus\upset\,(Q\cap M_\alpha))$.  Then (\ref{enumDstable}) is easily verified: if $q\in Q\cap M_\alpha$, then 
\begin{equation*}
D_{\alpha+1}\cap\upset q=(D_\alpha\cap\upset q)\cup((E\cap\upset q)\setminus\upset\,(Q\cap M_\alpha))=D_\alpha\cap\upset q.
\end{equation*}

Let us verify (\ref{enumDdense}) for stage $\alpha+1$.  Let $q\in Q\cap M_{\alpha+1}$. If $q\in\upset\,(Q\cap M_\alpha)$, then $q\in\upset D_\alpha\subseteq\upset D_{\alpha+1}$ because of (\ref{enumDdense}) for stage $\alpha$.  Suppose $q\not\in\upset\,(Q\cap M_\alpha)$.  Choose $e\in E$ such that $e\leq q$.  Then $e\not\in\upset\,(Q\cap M_\alpha)$; hence, $q\in\upset\,(E\setminus\upset\,(Q\cap M_\alpha))\subseteq\upset D_{\alpha+1}$.

It remains only to verify (\ref{enumDomegaop}) for stage $\alpha+1$.  Let $q\in Q\cap M_{\alpha+1}$.  Then $E\cap\upset q$ is finite; hence, by the definition of $D_{\alpha+1}$, it suffices to show that $D_\alpha\cap\upset q$ is finite.  By Lemma~\ref{LEMfreeboolreflect}, there exists $r\in B\cap M_\alpha$ such that $r\geq q$ and $M_\alpha\cap\upset q=M_\alpha\cap\upset r$; hence, $D_\alpha\cap\upset q=D_\alpha\cap\upset r$.  Since $q\in Q$, we have $r\in M_\alpha\cap\upset Q$.  By elementarity, there exists $p\in Q\cap M_\alpha$ such that $p\leq r$; hence, $D_\alpha\cap\upset r\subseteq D_\alpha\cap\upset p$.  By (\ref{enumDdense}) for stage $\alpha$, we have $D_\alpha\cap\upset p$ is finite; hence, $D_\alpha\cap\upset q$ is finite.
\end{proof}

\begin{definition}
For any space $X$, let $\Clop(X)$ denote the boolean algebra of clopen subsets of $X$.
\end{definition}

\begin{theorem}\label{THMalmostomegaopsubsetsdyadic}
Let $X$ be a dyadic compactum and let $\mcU$ be a family of subsets of $X$ such that for all $U\in\mcU$ there exists $V\in\mcU$ such that $\closure{V}\cap\closure{X\setminus U}=\emptyset$.  Then $\mcU$ is almost \omegaop.
\end{theorem}
\begin{proof}
Let $f\colon 2^\kappa\to X$ be a continuous surjection for some cardinal $\kappa$.  Set $\mcB=\Clop(2^\kappa)$.  Then $\mcB$ is a free boolean algebra.  Set $\mcV=\{\inv{f}U: U\in\mcU\}$.  Then it suffices to show that $\mcV$ is almost \omegaop.  Let $\mcQ$ denote the set of all $B\in\mcB$ such that $V\subseteq B$ for some $V\in\mcV$.  By Theorem~\ref{THMomegaopbool}, $\mcQ$ is almost \omegaop.  Hence, by Lemma~\ref{LEMmutuallydense}, it suffices to show that $\mcQ$ and $\mcV$ are mutually dense.  By definition, every $Q\in\mcQ$ contains some $V\in\mcV$; hence, it suffices to show that every $V\in\mcV$ contains some $Q\in\mcQ$.  Suppose $V\in\mcV$. Choose $U\in\mcU$ such that $\closure{U}\cap\closure{X\setminus f[V]}=\emptyset$.  Then there exists $B\in\mcB$ such that $\inv{f}\closure{U}\subseteq B\subseteq V$; hence, $V\supseteq B\in\mcQ$.
\end{proof}

The following corollary is immediate and it implies the first half of Theorem~\ref{THMdyadicmain}.

\begin{corollary}\label{CORdyadicpibaselocbase}
Let $X$ be a dyadic compactum.  Then, for all closed subsets $C$ of $X$, every neighborhood base of $C$ contains an \omegaop\ neighborhood base of $C$.  Moreover, every $\pi$\nbd-base of $X$ contains an \omegaop\ $\pi$\nbd-base of $X$.
\end{corollary}

\begin{remark}
The first half of the above corollary can also proved simply by citing Theorem~\ref{THMokcharctsprodw} and Lemma~\ref{LEMmutuallydense}.
\end{remark}

Next we state the natural generalizations of Lemma~\ref{LEMfreeboolreflect}, Theorem~\ref{THMomegaopbool}, Theorem~\ref{THMalmostomegaopsubsetsdyadic}, and Corollary~\ref{CORdyadicpibaselocbase} to continuous images of products of compacta with bounded weight.  We will only remark briefly about the proofs of these generalizations, for they are easy modifications of the corresponding old proofs.

\begin{lemma}\label{LEMcoprodboolreflect}
Let $\kappa$ be a regular uncountable cardinal and let $B$ be a coproduct $\coprod_{i\in I}B_i$ of boolean algebras all of size less than $\kappa$; let $\{B,\land,\lor,\la B_i\ra_{i\in I}\}\subseteq M\elemsub H_\theta$ and $M\cap\kappa\in\kappa+1$.  Then, for all $q\in B$, there exists $r\in B\cap M$ such that, for all $p\in B\cap M$, we have $p\geq q$ if and only if $p\geq r$.  In particular, $r\geq q$.
\end{lemma}
\begin{proof}
Note that the subalgebra $B\cap M$ is the subcoproduct $\coprod_{i\in I\cap M}B_i$ naturally embedded in $B$.  Then proceed as in the proof of Lemma~\ref{LEMfreeboolreflect} with $\bigcup_{i\in I}B_i$, naturally embedded in $B$, playing the role of $G$.
\end{proof}

\begin{theorem}\label{THMkappaopbool}
Let $\kappa\geq\omega$ and $B$ be a coproduct of boolean algebras all of size at most $\kappa$.  Then every subset of $B$ is almost \kappaop.
\end{theorem}
\begin{proof}
The proof is essentially the proof of Theorem~\ref{THMomegaopbool}.  Instead of using Lemma~\ref{LEMfreeboolreflect}, use the instance of Lemma~\ref{LEMcoprodboolreflect} for the regular uncountable cardinal $\kappa^+$.
\end{proof}

\begin{theorem}\label{THMalmostomegaopsubsetsgendyadic}
Let $\kappa\geq\omega$ and let $X$ be Hausdorff and a continuous image of a product of compacta all of weight at most $\kappa$; let $\mcU$ be a family of subsets of $X$ such that, for all $U\in\mcU$, there exists $V\in\mcU$ such that $\closure{V}\cap\closure{X\setminus U}=\emptyset$.  Then $\mcU$ is almost \kappaop.
\end{theorem}
\begin{proof}
Let $h\colon\prod_{i\in I}X_i\rightarrow X$ be a continuous surjection where each $X_i$ is a compactum with weight at most $\kappa$.  Each $X_i$ embeds into $[0,1]^\kappa$ and is therefore a continuous image of a closed subspace of $2^\kappa$.  Hence, \wma\ $\prod_{i\in I}X_i$ is totally disconnected.  The rest of the proof is just the proof of Theorem~\ref{THMalmostomegaopsubsetsdyadic} with Theorem~\ref{THMkappaopbool} replacing Theorem~\ref{THMomegaopbool}.
\end{proof}

The following corollary is immediate.

\begin{corollary}\label{CORgendyadicpibaselocbase}
Let $\kappa\geq\omega$ and let $X$ be Hausdorff and a continuous image of a product of compacta all of weight at most $\kappa$.  Then, for all closed subsets $C$ of $X$, every neighborhood base of $C$ contains a \kappaop\ neighborhood base of $C$.  Moreover, every $\pi$\nbd-base of $X$ contains a \kappaop\ $\pi$\nbd-base of $X$.
\end{corollary}

\begin{remark}
Again, the first half of the above corollary can also proved simply by citing Theorem~\ref{THMokcharctsprodw} and Lemma~\ref{LEMmutuallydense}.
\end{remark}

In contrast to Corollary~\ref{CORdyadicpibaselocbase}, not all dyadic compacta have \omegaop\ bases.  The following proposition is essentially due to Peregudov (see Lemma 1 of \cite{peregudov97}).  It makes it easy to produces examples of dyadic compacta $X$ such that $\ow{X}>\omega$.

\begin{proposition}\label{PROpicharcfchar}
Suppose a point $p$ in a space $X$ satisfies $\pi\character{p,X}<\cf\kappa=\kappa\leq\character{p,X}$.  Then $\ow{X}>\kappa$.
\end{proposition}
\begin{proof}
Let $\mcA$ be a base of $X$.  Let $\mcU_0$ and $\mcV_0$ be, respectively, a local $\pi$\nbd-base at $p$ of size at most $\pi\character{p,X}$ and a local base at $p$ of size $\character{p,X}$.  For each element of $\mcU_0$, choose a subset in $\mcA$, thereby producing a local $\pi$\nbd-base $\mcU$ at $p$ that is a subset of $\mcA$ of size at most $\pi\character{p,X}$.  Similarly, for each element of $\mcV_0$, choose a smaller neighborhood of $p$ in $\mcA$, thereby producing a local base $\mcV$ at $p$ that is a subset of $\mcA$ of size $\character{p,X}$.  Every element of $\mcV$ contains an element of $\mcU$.  Hence, some element of $\mcU$ is contained in $\kappa$\nbd-many elements of $\mcV$; hence, $\mcA$ is not \kappaop.
\end{proof}

\begin{example}\label{EXcharnotpichar}
Let $X$ be the discrete sum of $2^\omega$ and $2^{\omega_1}$.  Let $Y$ be the quotient of $X$ resulting from collapsing a point in $2^\omega$ and a point in $2^{\omega_1}$ to a single point $p$.  Then $\pi\character{p,Y}=\omega$ and $\character{p,Y}=\omega_1$; hence, $\ow{Y}>\omega_1$.
\end{example}

\begin{question}
Is there a dyadic compactum $X$ such that $\pi\character{p,X}=\character{p,X}$ for all $p\in X$ but $X$ has no \omegaop\ base?  In particular, if $Y$ is as in Example~\ref{EXcharnotpichar} and $Z$ is the discrete sum of $Y$ and $2^{\omega_2}$, then does $Z^{\omega_1}$ have an \omegaop\ base?
\end{question}

If we make an additional assumption about a dyadic compactum $X$, namely, that all its points have $\pi$\nbd-character equal to its weight, then $X$ has an \omegaop\ base.  Also, we may choose this \omegaop\ base to be a subset of an arbitrary base of $X$.  To prove this, we approximate such an $X$ by metric compacta.  Each such metric compactum is constructed using the following technique due to Bandlow~\cite{bandlow}.

\begin{definition}
Given a space $X$, let $C(X)$ denote the set of continuous maps from $X$ to $\reals$.
\end{definition}

\begin{definition}\label{DEFsubmodelquotient}
Suppose $X$ is a space and $\mcF\subseteq C(X)$.  For all $p\in X$, let $p/\mcF$ denote the set of $q\in X$ satisfying $f(p)=f(q)$ for all $f\in\mcF$.  For each $f\in\mcF$, define $f/\mcF\colon X/\mcF\to\reals$ by $(f/\mcF)(p/\mcF)=f(p)$ for all $p\in X$.
\end{definition}

\begin{lemma}\label{LEMquotienttop}
Suppose $X$ is a compactum and $\mcF\subseteq C(X)$.  Then $X/\mcF$ (with the quotient topology) is a compactum and its topology is the coarsest topology for which $f/\mcF$ is continuous for all $f\in\mcF$.  Further suppose $\{X\setminus\inv{f}\{0\}:f\in\mcF\}$ is a base of $X$ and $\mcF\in M\elemsub H_\theta$.  Then $\{(X\setminus\inv{f}\{0\})/(\mcF\cap M):f\in\mcF\cap M\}$ is a base of $X/(\mcF\cap M)$.
\end{lemma}
\begin{proof}
If $f\in\mcF$, then $f/\mcF$ is clearly continuous with respect to the quotient topology of $X/\mcF$.  Therefore, the compact quotient topology on $X/\mcF$ is finer than the Hausdorff topology induced by $\{f/\mcF: f\in\mcF\}$.  If a compact topology $\mcT_0$ is finer than a Hausdorff topology $\mcT_1$, then $\mcT_0=\mcT_1$.  Hence, the quotient topology on $X/\mcF$ is the topology induced by $\{f/\mcF: f\in\mcF\}$.  

Set $\mcA=\{X\setminus\inv{f}\{0\}:f\in\mcF\}$.  Suppose $\mcA$ is a base of $X$ and $\mcF\in M\elemsub H_\theta$.  Let us show that $\{(X\setminus\inv{f}\{0\})/(\mcF\cap M):f\in\mcF\cap M\}$ is a base of $X/(\mcF\cap M)$.  Let $\mcU$ denote the set of preimages of open rational intervals with respect to elements of $\mcF\cap M$.  Let $\mcV$ denote the set of nonempty finite intersections of elements of $\mcU$.  Then $\mcV\subseteq M$ and $\{V/(\mcF\cap M): V\in\mcV\}$ is base of $X/(\mcF\cap M)$.  Suppose $p\in V_0\in\mcV$.  Then it suffices to find $W\in\mcA\cap M$ such that $p\in W\subseteq V_0$.  Choose $V_1\in\mcV$ such that $p\in V_1\subseteq\closure{V}_1\subseteq V_0$.  Then there exist $n<\omega$ and $W_0,\ldots,W_{n-1}\in\mcA$ such that $\closure{V}_1\subseteq\bigcup_{i<n}W_i\subseteq V_0$.  By elementarity, \wma\ $W_0,\ldots,W_{n-1}\in M$.  Hence, there exists $i<n$ such that $p\in W_i\subseteq V_0$ and $W_i\in\mcA\cap M$.
\end{proof}

To construct an \omegaop\ base of a suitable dyadic compactum $X$, we apply Lemma~\ref{LEMcpctmetricbase} to a family of spaces $X/(\mcF\cap M)$ where $\mcF\subseteq C(X)$ and $M$ ranges over a transfinite sequence of countable elementary submodels of $H_\theta$.  This sequence is constructed such that, loosely speaking, each submodel in the sequence knows about the preceding submodels.

\begin{definition}\label{DEFapprox}
Let $\kappa$ be a regular uncountable cardinal and let $\la H_\theta,\ldots\ra$ be an expansion of the $\{\in\}$\nbd-structure $H_\theta$ to an $\mcL$\nbd-structure for some language $\mcL$ of size less than $\kappa$.  Then a $\kappa$\nbd-\emph{approximation sequence} in $\la H_\theta,\ldots\ra$ is an ordinally indexed sequence $\la M_\alpha\ra_{\alpha<\eta}$ such that for all $\alpha<\eta$ we have $\{\kappa,\la M_\beta\ra_{\beta<\alpha}\}\subseteq M_\alpha\elemsub\la H_\theta,\ldots\ra$ and $\card{M_\alpha}\subseteq M_\alpha\cap\kappa\in\kappa$.
\end{definition}

The following lemma is a generalization of a technique of Jackson and Mauldin~\cite{jackson} of approximating a structure by a tree of elementary substructures.

\begin{lemma}\label{LEMapprox}
Let $\kappa$ and $\la H_\theta,\ldots\ra$ be as in Definition~\ref{DEFapprox}.  Then there is a $\{\kappa\}$\nbd-definable map $\Psi$ that sends every $\kappa$\nbd-approximation sequence $\la M_\alpha\ra_{\alpha<\eta}$ in $\la H_\theta,\ldots\ra$ to a sequence $\la\Sigma_\alpha\ra_{\alpha\leq\eta}$ such that we have the following for all $\alpha\leq\eta$.
\begin{enumerate}
\item\label{enumfinitesigma} $\Sigma_\alpha$ is a finite set.
\item\label{enumsigmaelemsub} $\card{N}\subseteq N\elemsub\la H_\theta,\ldots\ra$ for all $N\in\Sigma_\alpha$.
\item\label{enumunionsigma} $\bigcup\Sigma_\alpha=\bigcup_{\beta<\alpha}M_\beta$.
\item\label{enumsigmainm} If $\alpha<\eta$, then $\Sigma_\alpha\in M_\alpha$.
\item\label{enumepsilonchain} $\Sigma_\alpha$ is an $\in$\nbd-chain.
\item\label{enumdesccard} If $N_0,N_1\in\Sigma_\alpha$ and $N_0\in N_1$, then $\card{N_0}>\card{N_1}$.
\item\label{enumpsisubseq} $\la\Sigma_\beta\ra_{\beta\leq\alpha}=\Psi(\la M_\beta\ra_{\beta<\alpha})$.
\end{enumerate}
Moreover, $\card{\Sigma_\lambda}=1$ and $\{\alpha<\lambda:\card{\Sigma_\alpha}=1\}$ is closed unbounded in $\lambda$ for all infinite cardinals $\lambda\leq\eta$.  
\end{lemma}
\begin{proof}
Let $\Omega$ denote the class of $\la\gamma_i\ra_{i<n}\in\ord^{<\omega}\setminus\{\emptyset\}$ for which $\kappa\leq\card{\gamma_i}>\card{\gamma_j}$ for all $i<j<n$ and $\card{\gamma_{n-1}}<\kappa$ if $n>0$.  Order $\Omega$ lexicographically and let $\Upsilon$ be the order isomorphism from $\ord$ to $\Omega$.  Given any $\sigma=\la\gamma_i\ra_{i<n}\in\ord^{<\omega}$ and $i<n$, set $\phi_i(\sigma)=\la\gamma_0,\ldots,\gamma_{i-1},0\ra$ and $\phi_n(\sigma)=\sigma$.  Let $\la M_\alpha\ra_{\alpha<\eta}$ be a $\kappa$\nbd-approximation sequence in $\la H_\theta,\ldots\ra$.  For all $\alpha\leq\eta$ and $i\in\dom\Upsilon(\alpha)$, set
\begin{equation*}
N_{\alpha,i}=\bigcup\{M_\beta:\phi_i(\Upsilon(\alpha))\leq\Upsilon(\beta)<\phi_{i+1}(\Upsilon(\alpha))\};
\end{equation*}
set $\Sigma_\alpha=\{N_{\alpha,i}:i\in\dom\Upsilon(\alpha)\}\setminus\{\emptyset\}$.  Then $\Psi$ is $\{\kappa\}$\nbd-definable and it is easily verified that $\card{\Sigma_\lambda}=1$ and $\{\alpha<\lambda:\card{\Sigma_\alpha}=1\}$ is closed unbounded in $\lambda$ for all infinite cardinals $\lambda\leq\eta$.  Let us prove (1)\nbd-(7).  (\ref{enumfinitesigma}), (\ref{enumunionsigma}), (\ref{enumsigmainm}), and (\ref{enumpsisubseq}) immediately follow from the relevant definitions.  Let $\alpha\leq\eta$ and $\la\beta_i\ra_{i<n}=\Upsilon(\alpha)$.  For all $\sigma\in\Omega$ and $i<n-1$, we have $\phi_i(\Upsilon(\alpha))\leq\sigma<\phi_{i+1}(\Upsilon(\alpha))$ if and only if $\sigma$ is the concatenation of $\la\beta_j\ra_{j<i}$ and some $\tau\in\Omega$ satisfying $\tau<\la\beta_i,0\ra$.  Therefore, $\card{N_{\alpha,i}}=\card{\beta_i}$ for all $i<n-1$.  
For all $\sigma\in\Omega$, we have $\phi_{n-1}(\Upsilon(\alpha))\leq\sigma<\phi_n(\Upsilon(\alpha))$ if and only if $\sigma=\la\beta_0,\ldots,\beta_{n-2},\gamma\ra$ for some $\gamma<\beta_{n-1}$.  Hence, $\card{N_{\alpha,n-1}}<\kappa$; hence, $\card{N_{\alpha,i}}>\card{N_{\alpha,j}}$ for all $i<j<n$.  Let $\Upsilon(\alpha_i)=\phi_i(\Upsilon(\alpha))$ for all $i<n$.  If $i<j<n$, then $\{N_{\alpha,k}:k<j\}=\Sigma_{\alpha_{j-1}}$; whence, either $N_{\alpha,j}=\emptyset$ or $N_{\alpha,i}\in M_{\alpha_{j-1}}\subseteq N_{\alpha,j}$, depending on whether $\beta_j=0$.  Thus, (\ref{enumepsilonchain}) and (\ref{enumdesccard}) hold.

Finally, let us prove (\ref{enumsigmaelemsub}).  Proceed by induction on $\alpha$.  Suppose $\beta_{n-1}>0$.  Since $\{N_{\alpha,i}:i<n-1\}=\Sigma_{\alpha_{n-1}}$ and $\alpha_{n-1}+\beta_{n-1}=\alpha$, it suffices to show that $\card{N_{\alpha,n-1}}\subseteq N_{\alpha,n-1}\elemsub\la H_\theta,\ldots\ra$.  If $\beta_{n-1}\in\Lim$, then $N_{\alpha,n-1}=\bigcup_{\gamma<\beta_{n-1}}N_{\alpha_{n-1}+\gamma,n-1}$; hence, $\card{N_{\alpha,n-1}}\subseteq N_{\alpha,n-1}\elemsub\la H_\theta,\ldots\ra$.  If $\beta_{n-1}\not\in\Lim$, then $N_{\alpha,n-1}=N_{\alpha-1,n-1}\cup M_{\alpha-1}=M_{\alpha-1}$ because $N_{\alpha-1,n-1}\in M_{\alpha-1}$ and $\card{N_{\alpha-1,n-1}}<\kappa$; hence, $\card{N_{\alpha,n-1}}\subseteq N_{\alpha,n-1}\elemsub\la H_\theta,\ldots\ra$.

Therefore, \wma\ $\beta_{n-1}=0$.  Hence, $\Sigma_\alpha=\{N_{\alpha,i}:i<n-1\}$; hence, \wma\ $n>1$.  Since $\{N_{\alpha,i}:i<n-2\}=\Sigma_{\alpha_{n-2}}$ and $\alpha_{n-2}<\alpha$, it suffices to show that $\card{N_{\alpha,n-2}}\subseteq N_{\alpha,n-2}\elemsub\la H_\theta,\ldots\ra$.  If $\beta_{n-2}=\kappa$, then $N_{\alpha,n-2}=\bigcup_{\gamma<\kappa}N_{\alpha_{n-2}+\gamma,n-2}$; hence, $\card{N_{\alpha,n-2}}\subseteq N_{\alpha,n-2}\elemsub\la H_\theta,\ldots\ra$.  Hence, \wma\ $\beta_{n-2}>\kappa$.  Let $\Upsilon(\delta_\gamma)=\la\beta_0,\ldots,\beta_{n-3},\gamma,0\ra$ for all $\gamma\in[\kappa,\beta_{n-2})$.  If $\beta_{n-2}\in\Lim$, then $N_{\alpha,n-2}=\bigcup_{\kappa\leq\gamma<\beta_{n-2}}N_{\delta_\gamma,n-2}$; hence, $\card{N_{\alpha,n-2}}\subseteq N_{\alpha,n-2}\elemsub\la H_\theta,\ldots\ra$.  Hence, we may let $\beta_{n-2}=\varepsilon+1$.  Suppose $\card{\varepsilon}=\kappa$.  Then $N_{\alpha,n-2}=N_{\delta_\varepsilon,n-2}\cup\bigcup_{\gamma<\kappa}M_{\delta_\varepsilon+\gamma}$.  If $\gamma<\kappa$, then $\phi_{n-1}(\Upsilon(\delta_\varepsilon+\gamma))=\Upsilon(\delta_\varepsilon)$; whence, $\delta_\varepsilon$ and $\gamma$ are definable from $\delta_\varepsilon+\gamma$ and $\kappa$; whence, $\gamma\cup\bigcup_{\rho<\gamma}M_{\delta_\varepsilon+\rho}\subseteq M_{\delta_\varepsilon+\gamma}$.  Hence, $\card{N_{\delta_\varepsilon,n-2}}=\kappa\subseteq\bigcup_{\gamma<\kappa}M_{\delta_\varepsilon+\gamma}\elemsub\la H_\theta,\ldots\ra$.  Moreover, since $N_{\delta_\varepsilon,n-2}\in M_{\delta_\varepsilon}$, we have $N_{\delta_\varepsilon,n-2}\subseteq\bigcup_{\gamma<\kappa}M_{\delta_\varepsilon+\gamma}$; hence, $\card{N_{\alpha,n-2}}=\kappa\subseteq N_{\alpha,n-2}\elemsub\la H_\theta,\ldots\ra$.  

Therefore, \wma\ $\card{\varepsilon}>\kappa$.  Let $\Upsilon(\zeta_\gamma)=\la\beta_0,\ldots,\beta_{n-3},\varepsilon,\kappa+\gamma,0\ra$ for all $\gamma<\card{\varepsilon}$.  Then $N_{\alpha,n-2}=N_{\delta_\varepsilon,n-2}\cup\bigcup_{\gamma<\card{\varepsilon}}N_{\zeta_\gamma,n-1}$.  If $\gamma<\card{\varepsilon}$, then $\Upsilon(\zeta_\gamma)(n-1)=\kappa+\gamma$; whence, $\gamma\in M_{\zeta_\gamma}\subseteq N_{\zeta_{\gamma+1},n-1}$.  Hence, $\card{\varepsilon}\subseteq\bigcup_{\gamma<\card{\varepsilon}}N_{\zeta_\gamma,n-1}\elemsub\la H_\theta,\ldots\ra$.  Since $\card{N_{\delta_\varepsilon,n-2}}=\card{\varepsilon}$ and $N_{\delta_\varepsilon,n-2}\in M_{\delta_\varepsilon}\subseteq N_{\zeta_0,n-1}$, we have $N_{\delta_\varepsilon,n-2}\subseteq\bigcup_{\gamma<\card{\varepsilon}}N_{\zeta_\gamma,n-1}$.  Hence, $\card{N_{\alpha,n-2}}=\card{\varepsilon}\subseteq N_{\alpha,n-2}\elemsub\la H_\theta,\ldots\ra$.
\end{proof}

\begin{proposition}\label{PROwnet}
If $X$ has a network consisting of at most $\weight{X}$\nbd-many closed subsets (in particular, if $X$ is regular), then every base of $X$ contains a base of size at most $\weight{X}$.
\end{proposition}
\begin{proof}
Let $\mcA$ be an arbitrary base of $X$; let $\mcB$ be a base of $X$ of size at most $\weight{X}$; let $\mcN$ a network of $X$ consisting of at most $\weight{X}$\nbd-many closed subsets.  Since $X$ is $\weight{X}^+$\nbd-compact, we may choose, for each $\la N,B\ra\in\mcN\times\mcB$ satisfying $N\subseteq B$, some $\mcU_{N,B}\in[\mcA]^{\leq\weight{X}}$ such that $N\subseteq\bigcup\mcU_{N,B}\subseteq B$.  Then $\bigcup\{\mcU_{N,B}:\mcN\ni N\subseteq B\in\mcB\}$ is a base of $X$ and in $[\mcA]^{\leq\weight{X}}$.
\end{proof}

\begin{lemma}\label{LEMweightpicharcozero}
Let $X$ be a dyadic compactum such that $\pi\character{p,X}=\weight{X}$ for all $p\in X$.  Let $\mcA$ be a base of $X$ consisting only of cozero sets.  Then $\mcA$ contains an \omegaop\ base of $X$.
\end{lemma}
\begin{proof}
Set $\kappa=\weight{X}$; by Proposition~\ref{PROwnet}, \wma\ $\card{\mcA}=\kappa$.  Choose $\mcF\subseteq C(X)$ such that $\mcA=\{X\setminus\inv{g}\{0\}:g\in\mcF\}$.  Let $h\colon 2^\lambda\to X$ be a continuous surjection for some cardinal $\lambda$.  Let $\mcB$ be the free boolean algebra $\Clop(2^\lambda)$.  By Lemma~\ref{LEMcpctmetricbase}, \wma\ $\kappa>\omega$.  Let $\la M_\alpha\ra_{\alpha<\kappa}$ be an $\omega_1$\nbd-approximation sequence in $\la H_\theta,\in,\mcF,h\ra$; set $\la\Sigma_\alpha\ra_{\alpha\leq\kappa}=\Psi(\la M_\alpha\ra_{\alpha<\kappa})$ as defined in Lemma~\ref{LEMapprox}.

For each $\alpha<\kappa$, set $\mcA_\alpha=\mcA\cap M_\alpha$ and $\mcF_\alpha=\mcF\cap M_\alpha$.  For every $\mcH\subseteq\mcA_\alpha$, let $\mcH/\mcF_\alpha$ denote $\{U/\mcF_\alpha:U\in\mcH\}$.  By Lemma~\ref{LEMquotienttop}, $\mcA_\alpha/\mcF_\alpha$ is a base of $X/\mcF_\alpha$.  Since $X/\mcF_\alpha$ is a metric compactum, there exists $\mcW_\alpha\subseteq\mcA_\alpha$ such that $\mcW_\alpha/\mcF_\alpha$ is a base of $X/\mcF_\alpha$ satisfying (\ref{enumMetricomegaop}), (\ref{enumMetricclose}), and (\ref{enumMetricfin}) of Lemma~\ref{LEMcpctmetricbase}.  By (\ref{enumMetricomegaop}) of Lemma~\ref{LEMcpctmetricbase}, we may choose, for each $U\in\mcW_\alpha$, some $E_{\alpha,U}\in\mcB\cap M_\alpha$ such that $\inv{h}\closure{U}\subseteq E_{\alpha,U}\subseteq\inv{h}V$ for all $V\in\mcW_\alpha$ satisfying $\closure{U}\subseteq V$.  Set $\mcG_\alpha=\{E_{\alpha,U}:U\in\mcW_\alpha\}$.

Suppose $\mcG_\alpha$ is not \omegaop.  Then there exist $U\in\mcW_\alpha$ and $\la V_n\ra_{n<\omega}\in\mcW_\alpha^\omega$ such that $E_{\alpha,U}\subsetneq E_{\alpha,V_n}\not=E_{\alpha,V_m}$ for all $m<n<\omega$.  Set $\Gamma=\{W\in\mcW_\alpha:U\subsetneq W\}$.  By (\ref{enumMetricomegaop}) of Lemma~\ref{LEMcpctmetricbase}, $\Gamma$ is finite; hence, by (\ref{enumMetricfin}) of Lemma~\ref{LEMcpctmetricbase}, there exists $n<\omega$ such that $\{W\in\mcW_\alpha:V_n\subsetneq W\}\not\subseteq\Gamma$.  Hence, there exists $W\in\mcW_\alpha$ such that $W$ strictly contains $V_n$ but not $U$.  Hence, by (\ref{enumMetricclose}) of Lemma~\ref{LEMcpctmetricbase}, $E_{\alpha,V_n}\subseteq\inv{h}W$; hence, $\inv{h}U\subseteq E_{\alpha,U}\subsetneq E_{\alpha,V_n}\subseteq\inv{h}W$; hence, $U\subsetneq W$, which is absurd.  Therefore, $\mcG_\alpha$ is \omegaop.

Let $\mcV_\alpha$ denote the set of $V\in\mcW_\alpha$ satisfying $U\not\subseteq V$ for all nonempty open $U\in\bigcup\Sigma_\alpha$.  Let us show that $\mcV_\alpha/\mcF_\alpha$ is a base of $X/\mcF_\alpha$.  If $V\in\mcV_\alpha$, then $\powset{V}\cap\mcW_\alpha\subseteq\mcV_\alpha$; hence, it suffices to show that $\mcV_\alpha$ covers $X$.  Since $\card{\bigcup\Sigma_\alpha}<\kappa$, every point of $X$ has a neighborhood in $\mcA$ that does not contain any nonempty open subset of $X$ in $\bigcup\Sigma_\alpha$.  By compactness, there is cover of $X$ by finitely many such neighborhoods, say, $W_0,\ldots,W_{n-1}$.  By elementarity, \wma\ $W_0,\ldots,W_{n-1}\in\mcA_\alpha$.  Then $\{W_i:i<n\}$ has a refining cover $\mcS\subseteq\mcW_\alpha$.  Hence, $\mcS\subseteq\mcV_\alpha$; hence, $\mcV_\alpha$ covers $X$ as desired.

Let $\mcU_\alpha$ denote the set of $U\in\mcV_\alpha$ such that $\closure{U}\subseteq V$ for some $V\in\mcV_\alpha$.  Then $\mcU_\alpha/\mcF_\alpha$ is clearly a base of $X/\mcF_\alpha$.  Set $\mcE_\alpha=\{E_{\alpha,U}:U\in\mcU_\alpha\}$.  Then $\mcE_\alpha$ is \omegaop\ because it is a subset of $\mcG_\alpha$.

For all $\mcI\subseteq\powset{2^\kappa}$, set $\upset\mcI=\{H\subseteq 2^\kappa:H\supseteq I\text{ for some }I\in\mcI\}$.  For all $H\subseteq 2^\kappa$, set $\upset H=\upset\,\{H\}$.  Set $\mcU=\bigcup_{\alpha<\kappa}\mcU_\alpha$ and $\mcC=\mcB\cap\upset\,\{\inv{h}U: U\in\mcU\}$.  For all $\alpha\leq\kappa$, set $\mcD_\alpha=\bigcup_{\beta<\alpha}\mcE_\beta$.  Then we claim the following for all $\alpha\leq\kappa$.
\begin{enumerate}
\item\label{enumDdense1} $\mcD_\alpha$ is a dense subset of $\mcC\cap\bigcup\Sigma_\alpha$.
\item\label{enumDomegaop1} $\mcD_\alpha\cap\upset H$ is finite for all $H\in\mcC\cap\bigcup\Sigma_\alpha$.
\item\label{enumDstable1} If $\alpha<\kappa$, then $\mcD_{\alpha+1}\cap\upset H=\mcD_\alpha\cap\upset H$ for all $H\in\mcC\cap\bigcup\Sigma_\alpha$.
\end{enumerate}
We prove this claim by induction.  For stage $0$, the claim is vacuous.  For limit stages, (\ref{enumDdense1}) is clearly preserved, and (\ref{enumDomegaop1}) is preserved because of (\ref{enumDstable1}).  Suppose $\alpha<\kappa$ and (\ref{enumDdense1}) and (\ref{enumDomegaop1}) hold for stage $\alpha$.  Then it suffices to prove (\ref{enumDstable1}) for stage $\alpha$ and to prove (\ref{enumDdense1}) and (\ref{enumDomegaop1}) for stage $\alpha+1$.  

Let us verify (\ref{enumDstable1}).  Seeking a contradiction, suppose $H\in\mcC\cap\bigcup\Sigma_\alpha$ and $\mcD_{\alpha+1}\cap\upset H\not=\mcD_\alpha\cap\upset H$.  Then $\mcE_\alpha\cap\upset H\not=\emptyset$; hence, there exists $U\in\mcU_\alpha$ such that $H\subseteq E_{\alpha,U}$.  By (\ref{enumDdense1}), there exist $\beta<\alpha$ and $W\in\mcU_\beta$ such that $E_{\beta,W}\subseteq H$.  By definition, there exists $V\in\mcV_\alpha$ such that $\closure{U}\subseteq V$.  Hence, $\inv{h}W\subseteq E_{\beta,W}\subseteq H\subseteq E_{\alpha,U}\subseteq\inv{h}V$; hence, $W\subseteq V$.  Since $W\in M_\beta\subseteq\bigcup\Sigma_\alpha$ and $V\in\mcV_\alpha$, we have $W\not\subseteq V$, which yields our desired contradiction.

Let us verify (\ref{enumDdense1}) for stage $\alpha+1$.  By (\ref{enumDdense1}) for stage $\alpha$, we have
\begin{equation*}
\mcD_{\alpha+1}=\mcD_\alpha\cup\mcE_\alpha\subseteq\Bigl(\mcC\cap\bigcup\Sigma_\alpha\Bigr)\cup(\mcC\cap M_\alpha)=\mcC\cap\bigcup\Sigma_{\alpha+1}, 
\end{equation*}
so we just need to show denseness.  Let $H\in\mcC\cap\bigcup\Sigma_{\alpha+1}$.  If $H\in\bigcup\Sigma_\alpha$, then $H\in\upset\mcD_\alpha$, so \wma\ $H\in M_\alpha$.  By elementarity, there exists $U_0\in\mcU_\alpha$ such that $\inv{h}U_0\subseteq H$.  Choose $U_1\in\mcU_\alpha$ such that $\closure{U}_1\subseteq U_0$.  Then $E_{\alpha,U_1}\subseteq\inv{h}U_0$; hence, $E_{\alpha,U_1}\subseteq H$.  Hence, $H\in\upset\mcD_{\alpha+1}$.

To complete the proof of the claim, let us verify (\ref{enumDomegaop1}) for stage $\alpha+1$.  By (\ref{enumDdense1}) for stage $\alpha+1$, it suffices to prove $\mcD_{\alpha+1}\cap\upset H$ is finite for all $H\in\mcD_{\alpha+1}$.  By (\ref{enumDstable1}), if $H\in\mcD_\alpha$, then $\mcD_{\alpha+1}\cap\upset H=\mcD_\alpha\cap\upset H$, which is finite by (\ref{enumDdense1}) and (\ref{enumDomegaop1}) for stage $\alpha$.  Hence, \wma\ $H\in\mcE_\alpha$.  Since $\mcE_\alpha$ is \omegaop, it suffices to show that $\mcD_\alpha\cap\upset H$ is finite.  Since $\mcD_\alpha\subseteq\bigcup\Sigma_\alpha$, it suffices to show that $\mcD_\alpha\cap N\cap\upset H$ is finite for all $N\in\Sigma_\alpha$.  Let $N\in\Sigma_\alpha$.  By Lemma~\ref{LEMfreeboolreflect}, there exists $G\in\mcB\cap N$ such that $G\supseteq H$ and $\mcB\cap N\cap\upset H=\mcB\cap N\cap\upset G$; hence, $\mcD_\alpha\cap N\cap\upset H=\mcD_\alpha\cap N\cap\upset G$.  Since $G\supseteq H\in\mcC$, we have $G\in\mcC$.  By (\ref{enumDomegaop1}) for stage $\alpha$, the set $\mcD_\alpha\cap N\cap\upset G$ is finite; hence, $\mcD_\alpha\cap N\cap\upset H$ is finite.  

Since $\mcU\subseteq\mcA$, it suffices to prove that $\mcU$ is an \omegaop\ base of $X$.  Suppose $p\in V\in\mcA$.  Then there exists $\alpha<\kappa$ such that $V\in\mcA_\alpha$.  Hence, there exists $U\in\mcU_\alpha$ such that $p/\mcF_\alpha\in U/\mcF_\alpha\subseteq V/\mcF_\alpha$; hence, $p\in U\subseteq V$.  Thus, $\mcU$ is a base of $X$.

Let us show that $\mcU$ is \omegaop.  Suppose not.  Then there exists $\alpha<\kappa$ and $U_0\in\mcU_\alpha$ such that there exist infinitely many $V\in\mcU$ such that $U_0\subseteq V$.  Choose $U_1\in\mcU_\alpha$ such that $\closure{U}_1\subseteq U_0$.  Suppose $\beta<\kappa$ and $U_0\subseteq V\in\mcU_\beta$.  Then $E_{\alpha,U_1}\subseteq\inv{h}U_0\subseteq\inv{h}V\subseteq E_{\beta,V}$.  By (\ref{enumDdense1}) and (\ref{enumDomegaop1}), $\mcD_\kappa$ is \omegaop; hence, there are only finitely many possible values for $E_{\beta,V}$.  Therefore, there exist $\la\gamma_n\ra_{n<\omega}\in\kappa^\omega$ and $\la V_n\ra_{n<\omega}\in\prod_{n<\omega}\mcU_{\gamma_n}$ such that $V_m\not=V_n$ and $E_{\gamma_m,V_m}=E_{\gamma_n,V_n}$ for all $m<n<\omega$.  Suppose that for some $\delta<\kappa$ we have $\gamma_n=\delta$ for all $n<\omega$.  Let $i<\omega$ and set $\Gamma=\{W\in\mcW_\delta: V_i\subsetneq W\}$.  By (\ref{enumMetricomegaop}) and (\ref{enumMetricfin}) of Lemma~\ref{LEMcpctmetricbase}, there exists $j<\omega$ such that $\{W\in\mcW_\delta: V_j\subsetneq W\}\not\subseteq\Gamma$.  Hence, there exists $W\in\mcW_\delta$ such that $W$ strictly contains $V_j$ but not $V_i$.  By (\ref{enumMetricclose}) of Lemma~\ref{LEMcpctmetricbase}, $\closure{V}_j\subseteq W$.  Hence, $\inv{h}\closure{V}_i\subseteq E_{\delta,V_i}=E_{\delta,V_j}\subseteq\inv{h}W$.  Hence, $\closure{V}_i\subseteq W$.  Since $W$ does not strictly contain $V_i$, we must have $V_i=\closure{V}_i=W$.  Hence, $\inv{h}V_i=E_{\delta,V_i}=E_{\delta,V_0}$.  Since $i$ was arbitrary chosen, we have $V_m=V_n=h[E_{\delta,V_0}]$ for all $m,n<\omega$, which is absurd.  Therefore, our supposed $\delta$ does not exist; hence, \wma\ $\gamma_0<\gamma_1$.  By definition, there exists $W\in\mcV_{\gamma_1}$ such that $\closure{V}_1\subseteq W$.  Therefore, $\inv{h}V_0\subseteq E_{\gamma_0,V_0}=E_{\gamma_1,V_1}\subseteq\inv{h}W$; hence, $V_0\subseteq W$.  Since $V_0\in M_{\gamma_0}\subseteq\bigcup\Sigma_{\gamma_1}$ and $W\in\mcV_{\gamma_1}$, we have $V_0\not\subseteq W$, which is absurd.  Therefore, $\mcU$ is \omegaop.
\end{proof}

Let us show that we may remove the requirement that the base $\mcA$ in Lemma~\ref{LEMweightpicharcozero} consist only of cozero sets.

\begin{proposition}\label{PROowcardbase}
If $X$ is a space and $\mcA$ is a \cardop{(\weight{X}^+)}\ base of $X$, then $\card{\mcA}\leq\weight{X}$.
\end{proposition}
\begin{proof}
Seeking a contradiction, suppose $\card{\mcA}>\weight{X}$.  Let $\mcB$ be a base of $X$ of size $\weight{X}$.  Then every element of $\mcA$ contains an element of $\mcB$.  Hence, some $U\in\mcB$ is contained in $\weight{X}^+$\nbd-many elements of $\mcA$.  Clearly $U$ contains some $V\in\mcA$, so $\mcA$ is not \cardop{(\weight{X}^+)}.
\end{proof}
\begin{remark}
The above proposition's proof can be trivially modified to show that if $\mcA$ is a \cardop{(\piweight{X}^+)}\ $\pi$\nbd-base of $X$, then $\card{\mcA}\leq\piweight{X}$.  Likewise, if $p\in X$ and $\mcA$ is a \cardop{(\character{p,X}^+)}\ local base at $p$, then $\card{\mcA}\leq\character{p,X}$.
\end{remark}

\begin{lemma}\label{LEMweightcharsubsetbase}
Suppose $X$ is a space with no isolated points and $\character{p,X}=\weight{X}$ for all $p\in X$.  Further suppose $\kappa=\cf\kappa\leq\min\{\ow{X},\weight{X}\}$ and $X$ has a network consisting of at most $\weight{X}$\nbd-many $\kappa$\nbd-compact sets.  Then every base of $X$ contains an \cardop{\ow{X}}\ base of $X$.
\end{lemma}
\begin{proof}
Set $\lambda=\ow{X}$ and $\mu=\weight{X}$.  Let $\mcA$ be an arbitrary base of $X$; let $\mcB$ be a \cardop{\lambda}\ base of $X$; let $\mcN$ be a network of $X$ consisting of at most $\mu$\nbd-many $\kappa$\nbd-compact sets.  By Proposition~\ref{PROowcardbase}, $\card{\mcB}\leq\mu$.  Let $\la\la N_\alpha,B_\alpha\ra\ra_{\alpha<\mu}$ enumerate $\{\la N,B\ra\in\mcN\times\mcB: N\subseteq B\}$.  Construct a sequence $\la\mcG_\alpha\ra_{\alpha<\mu}$ as follows.  Suppose $\alpha<\mu$ and $\la\mcG_\beta\ra_{\beta<\alpha}$ is a sequence of elements of $[\mcB]^{<\kappa}$.  For each $p\in N_\alpha$, we have $\character{p,X}=\mu\geq\kappa=\cf\kappa$; hence, we may choose $U_{\alpha,p}\in\mcB$ such that $p\in U_{\alpha,p}\not\in\bigcup_{\beta<\alpha}\mcG_\beta$.  Choose $\sigma_\alpha\in\bigl[\,N_\alpha\bigr]^{<\kappa}$ such that $N_\alpha\subseteq\bigcup_{p\in\sigma_\alpha}U_{\alpha,p}$.  Set $\mcG_\alpha=\{U_{\alpha,p}:p\in\sigma_\alpha\}$.

For each $\alpha<\mu$, choose $\mcF_\alpha\in[\mcA]^{<\kappa}$ such that $N_\alpha\subseteq\bigcup\mcF_\alpha\subseteq B_\alpha$ and $\mcF_\alpha$ refines $\mcG_\alpha$.  Set $\mcF=\bigcup_{\alpha<\mu}\mcF_\alpha$, which is easily seen to be a base of $X$.  Let us show that $\mcF$ is \cardop{\lambda}.  Suppose not.  Then, since $\kappa=\cf\kappa\leq\lambda$, there exist $V\in\mcF$ and $I\in[\mu]^\lambda$ and $\la W_\alpha\ra_{\alpha\in I}\in\prod_{\alpha\in I}\mcF_{\alpha}$ such that $V\subseteq\bigcap_{\alpha\in I}W_\alpha$.  For each $\alpha\in I$, there is a superset of $W_\alpha$ in $\mcG_\alpha$.  By induction, $\mcG_\alpha\cap\mcG_\beta=\emptyset$ for all $\alpha<\beta<\mu$; hence, $V$ has $\lambda$\nbd-many supersets in the \cardop{\lambda}\ base $\mcB$, which is absurd, for $V$ has a subset in $\mcB$.
\end{proof}

\begin{remark}
If $X$ is regular and locally $\kappa$\nbd-compact and $\kappa\leq\weight{X}$, then it is easily seen that $X$ has a network consisting of at most $\weight{X}$\nbd-many $\kappa$\nbd-compact closed sets.
\end{remark}

\begin{theorem}\label{THMweightpichar}
Let $X$ be a dyadic compactum such that $\pi\character{p,X}=\weight{X}$ for all $p\in X$.  Then every base $\mcA$ of $X$ contains an \omegaop\ base of $X$.
\end{theorem}
\begin{proof}
By Lemma~\ref{LEMweightpicharcozero}, $\ow{X}=\omega$.  Since $\weight{X}=\pi\character{p,X}\leq\character{p,X}\leq\weight{X}$ for all $p\in X$, we may apply Lemma~\ref{LEMweightcharsubsetbase} to get a subset of $\mcA$ that is an \omegaop\ base of $X$.
\end{proof}

Finally, let us prove the second half of Theorem~\ref{THMdyadicmain}.

\begin{corollary}\label{CORhomogdyadic}
Let $X$ be a homogeneous dyadic compactum with base $\mcA$.  Then $\mcA$ contains an \omegaop\ base of $X$.
\end{corollary}
\begin{proof}
Efimov~\cite{efimov} and Gerlits~\cite{gerlits} independently proved that the $\pi$\nbd-character of every dyadic compactum is equal to its weight.  Since $X$ is homogeneous, $\pi\character{p,X}=\weight{X}$ for all $p\in X$.  Hence, $\mcA$ contains an \omegaop\ base of $X$ by Theorem~\ref{THMweightpichar}.
\end{proof}

Note that a compactum is dyadic if and only if it a continous image of a product of second countable compacta.  Let us prove generalizations of Theorem~\ref{THMweightpichar} and Corollary~\ref{CORhomogdyadic} about continuous images of products of compacta with bounded weight.

\begin{lemma}\label{LEMweightequalpichar}
Suppose $\kappa=\cf\kappa>\omega$ and $X$ is a space such that $\pi\character{p,X}=\weight{X}\geq\kappa$ for all $p\in X$.  Further suppose $X$ has a network consisting of at most $\weight{X}$\nbd-many $\kappa$\nbd-compact closed sets.  Then every base of $X$ contains a \cardop{\weight{X}}\ base of $X$.
\end{lemma}
\begin{proof}
Set $\lambda=\weight{X}$ and let $\mcA$ be an arbitrary base of $X$.  By Proposition~\ref{PROwnet}, \wma\ $\card{\mcA}=\lambda$.  Let $\mcN$ be a network of $X$ consisting of at most $\lambda$\nbd-many $\kappa$\nbd-compact sets.  Let $\la M_\alpha\ra_{\alpha<\lambda}$ be a continuous elementary chain such that $\mcA,\mcN,M_\alpha\in M_{\alpha+1}\elemsub H_\theta$ and for all $\alpha<\lambda$.  We may also require that $M_{\alpha}\cap\kappa\in\kappa>\card{M_\alpha}$ for all $\alpha<\kappa$ and $\card{M_\alpha}=\card{\kappa+\alpha}$ for all $\alpha\in\lambda\setminus\kappa$.  For each $\alpha<\lambda$, set $\mcA_\alpha=\mcA\cap M_\alpha$.  Set $\mcB=\bigcup_{\alpha<\lambda}\mcA_{\alpha+1}\setminus\upset\mcA_\alpha$, which is clearly \cardop{\lambda}.  Let us show that $\mcB$ is a base of $X$.  Suppose $p\in U\in\mcA$.  Choose $N\in\mcN$ such that $p\in N\subseteq U$.  Choose $\alpha<\lambda$ such that $N,U\in\mcA_{\alpha+1}$.  For each $q\in N$, choose $V_q\in\mcA\setminus\upset\mcA_\alpha$ such that $q\in V_q\subseteq U$.  Then there exists $\sigma\in[N]^{<\kappa}$ such that $N\subseteq\bigcup_{q\in\sigma}V_q$.  By elementarity, \wma\ $\la V_q\ra_{q\in\sigma}\in M_{\alpha+1}$.  Choose $q\in\sigma$ such that $p\in V_q$.  Then $V_q\in\mcB$ and $p\in V_q\subseteq U$.  Thus, $\mcB$ is a base of $X$.
\end{proof}

\begin{theorem}\label{THMgendyadicweightpichar}
Let $\kappa\geq\omega$ and let $X$ be Hausdorff and a continuous image of a product of compacta each with weight at most $\kappa$.  Suppose $\pi\character{p,X}=\weight{X}$ for all $p\in X$.  Then every base of $X$ contains a \kappaop\ base.
\end{theorem}
\begin{proof}
Let $h\colon\prod_{i\in I}X_i\rightarrow X$ be a continuous surjection where each $X_i$ is a compactum with weight at most $\kappa$.  Each $X_i$ embeds into $[0,1]^\kappa$ and is therefore a continuous image of a closed subspace of $2^\kappa$.  Hence, \wma\ $\prod_{i\in I}X_i$ is totally disconnected.  Set $\lambda=\weight{X}$; by Lemmas~\ref{LEMcpctmetricbase} and \ref{LEMweightequalpichar}, \wma\ $\lambda>\kappa$.  By Theorem~\ref{THMweightpichar}, \wma\ $\kappa>\omega$.  Inductively construct a $\kappa^+$\nbd-approximation sequence $\la M_\alpha\ra_{\alpha<\lambda}$ in $\la H_\theta,\in,C(X),h,\la\Clop(X_i)\ra_{i\in I}\ra$ as follows.  For each $\alpha<\lambda$, let $\la N_{\alpha,\beta}\ra_{\beta<\kappa}$ be an $\omega_1$\nbd-approximation sequence in
\begin{equation*}
\la H_\theta,\in,C(X),h,\kappa,\la\Clop(X_i)\ra_{i\in I}, \la M_\beta\ra_{\beta<\alpha}\ra.
\end{equation*}
Set $\la\Gamma_{\alpha,\beta}\ra_{\beta\leq\kappa}=\Psi(\la N_{\alpha,\beta}\ra_{\beta<\kappa})$ as defined in Lemma~\ref{LEMapprox}; let $\{M_\alpha\}=\Gamma_{\alpha,\kappa}$.  Set $\la\Sigma_\alpha\ra_{\alpha\leq\lambda}=\Psi(\la M_\alpha\ra_{\alpha<\lambda})$.  Set $\mcF=C(X)\cap\bigcup\Sigma_\lambda$ and $\mcA=\{X\setminus\inv{f}\{0\}:f\in\mcF\}$.  Then $\mcA$ is a base of $X$.  By Lemma~\ref{LEMweightcharsubsetbase}, it suffices to construct a subset of $\mcA$ that is a \kappaop\ base of $X$.

For each $\alpha<\lambda$, set $\mcF_\alpha=\mcF\cap M_\alpha$.  Let $\mcV_\alpha$ denote the set of $V\in\mcA\cap M_\alpha$ satisfying $U\not\subseteq V$ for all nonempty open $U\in\bigcup\Sigma_\alpha$.  Arguing as in the proof Lemma~\ref{LEMweightpicharcozero}, $\mcV_\alpha/\mcF_\alpha$ is a base of $X/\mcF_\alpha$.  For each $\beta<\kappa$, let $\mcV_{\alpha,\beta}$ denote the set of all $V\in\mcV_\alpha\cap N_{\alpha,\beta}$ satisfying $U\not\subseteq V$ for all nonempty open $U\in\bigcup\Gamma_{\alpha,\beta}$.  Let $\mcR_{\alpha,\beta}$ denote the set of $\la U,V\ra\in\mcV_{\alpha,\beta}^2$ for which $\closure{U}\subseteq V$; set $\mcU_{\alpha,\beta}=\dom\mcR_{\alpha,\beta}$; set $\mcU_\alpha=\bigcup_{\beta<\kappa}\mcU_{\alpha,\beta}$.

Let us show that $\mcU_\alpha/\mcF_\alpha$ is also a base of $X/\mcF_\alpha$.  Suppose $p\in V\in\mcV_\alpha$.  Extend $\{V\}$ to a finite subcover $\sigma$ of $\mcV_\alpha$ such that $p\not\in\bigcup(\sigma\setminus\{V\})$.  Choose $\beta<\kappa$ such that $\sigma\in N_{\alpha,\beta}$.  For each $q\in X$, choose $V_{q,0},V_{q,1}\in\mcA$ such that $q\in V_{q,0}$ and there exists $W\in\sigma$ such that $U\not\subseteq\closure{V}_{q,0}\subseteq V_{q,1}\subseteq W$ for all nonempty open $U\in\bigcup\Sigma_\alpha\cup\bigcup\Gamma_{\alpha,\beta}$.  Choose $\tau\in[X]^{<\omega}$ such that $X=\bigcup_{q\in\tau}V_{q,0}$.  By elementarity, \wma\ $\la V_{q,i}\ra_{\la q,i\ra\in\tau\times 2}\in N_{\alpha,\beta}$.  Choose $q\in\tau$ such that $p\in V_{q,0}$.  Then $V_{q,0}\in\mcU_{\alpha,\beta}$ and $p\in V_{q,0}\subseteq V$.  Thus, $\mcU_\alpha/\mcF_\alpha$ is a base of $X/\mcF_\alpha$.

Set $\mcB=\Clop\bigl(\prod_{i\in I}X_i\bigr)$.  For each $\la U_0,U_1\ra\in\bigcup_{\beta<\kappa}\mcR_{\alpha,\beta}$, choose $E_\alpha(U_0,U_1)\in\mcB\cap M_\alpha$ such that $\inv{h}\closure{U}_0\subseteq E_\alpha(U_0,U_1)\subseteq\inv{h}U_1$.  Set $\mcE_{\alpha,\beta}=E_\alpha[\mcR_{\alpha,\beta}]$.  Set $\mcE_\alpha=\bigcup_{\beta<\kappa}\mcE_{\alpha,\beta}$.  Let us show that $\mcE_\alpha$ is \kappaop.  Suppose $\beta,\gamma<\kappa$ and $\mcE_{\alpha,\beta}\ni H\subseteq K\in\mcE_{\alpha,\gamma}$.  Then it suffices to show that $\gamma\leq\beta$.  Seeking a contradiction, suppose $\beta<\gamma$.  There exist $\la U_0,U_1\ra\in\mcR_{\alpha,\beta}$ and $\la V_0,V_1\ra\in\mcR_{\alpha,\gamma}$ such that $H=E_\alpha(U_0,U_1)$ and $K=E_\alpha(V_0,V_1)$.  Hence, $\bigcup\Gamma_{\alpha,\gamma}\ni U_0\subseteq V_1\in\mcV_{\alpha,\gamma}$, in contradiction with the definition of $\mcV_{\alpha,\gamma}$.

Set $\mcU=\bigcup_{\alpha<\lambda}\mcU_\alpha$ and $\mcC=\mcB\cap\upset\{\inv{h}U:U\in\mcU\}$.  For all $\alpha\leq\lambda$, set $\mcD_\alpha=\bigcup_{\beta<\alpha}\mcE_\beta$.  Then we claim the following for all $\alpha\leq\lambda$.
\begin{enumerate}
\item\label{enumDdense3} $\mcD_\alpha$ is a dense subset of $\mcC\cap\bigcup\Sigma_\alpha$.
\item\label{enumDkappaop3} $\card{\mcD_\alpha\cap\upset H}<\kappa$ for all $H\in\mcC\cap\bigcup\Sigma_\alpha$.
\item\label{enumDstable3} If $\alpha<\lambda$, then $\mcD_{\alpha+1}\cap\upset H=\mcD_\alpha\cap\upset H$ for all $H\in\mcC\cap\bigcup\Sigma_\alpha$.
\end{enumerate}
We prove this claim by induction.  For stage $0$, the claim is vacuous.  For limit stages, (\ref{enumDdense3}) is clearly preserved, and (\ref{enumDkappaop3}) is preserved because of (\ref{enumDstable3}).  Suppose $\alpha<\kappa$ and (\ref{enumDdense3}) and (\ref{enumDkappaop3}) hold for stage $\alpha$.  Then it suffices to prove (\ref{enumDstable3}) for stage $\alpha$ and to prove (\ref{enumDdense3}) and (\ref{enumDkappaop3}) for stage $\alpha+1$.  

Let us verify (\ref{enumDstable3}).  Seeking a contradiction, suppose $H\in\mcC\cap\bigcup\Sigma_\alpha$ and $\mcD_{\alpha+1}\cap\upset H\not=\mcD_\alpha\cap\upset H$.  Then $\mcE_\alpha\cap\upset H\not=\emptyset$; hence, there exists $V\in\mcU_\alpha$ such that $H\subseteq\inv{h}V$.  By (\ref{enumDdense3}), there exist $\beta<\alpha$ and $U\in\mcU_\beta$ and $K\in\mcE_\beta$ such that $\inv{h}\closure{U}\subseteq K\subseteq H$.  Hence, $U\subseteq V$.  Since $U\in M_\beta\subseteq\bigcup\Sigma_\alpha$ and $V\in\mcV_\alpha$, we have $U\not\subseteq V$, which yields our desired contradiction.

Let us verify (\ref{enumDdense3}) for stage $\alpha+1$.  By (\ref{enumDdense3}) for stage $\alpha$, we have
\begin{equation*}
\mcD_{\alpha+1}=\mcD_\alpha\cup\mcE_\alpha\subseteq\Bigl(\mcC\cap\bigcup\Sigma_\alpha\Bigr)\cup(\mcC\cap M_\alpha)=\mcC\cap\bigcup\Sigma_{\alpha+1}, 
\end{equation*}
so we just need to show denseness.  Let $H\in\mcC\cap\bigcup\Sigma_{\alpha+1}$.  If $H\in\bigcup\Sigma_\alpha$, then $H\in\upset\mcD_\alpha$, so \wma\ $H\in M_\alpha$.  By elementarity, there exists $U\in\mcU_\alpha$ such that $\inv{h}U\subseteq H$.  Choose $\beta<\kappa$ such that $U\in\mcU_{\alpha,\beta}$; choose $V\in\mcU_{\alpha,\beta}$ such that $\closure{V}\subseteq U$.  Then $E_\alpha(V,U)\subseteq H$; hence, $H\in\upset\mcD_{\alpha+1}$.

The proof of the claim is completed by noting that (\ref{enumDkappaop3}) for stage $\alpha+1$ can be verfied just as in the proof of Lemma~\ref{LEMweightpicharcozero}, except that Lemma~\ref{LEMcoprodboolreflect} is used in place of Lemma~\ref{LEMfreeboolreflect}.

Just as in the proof of Lemma~\ref{LEMweightpicharcozero}, $\mcU$ is a base of $X$; hence, it suffices to show that $\mcU$ is \kappaop.  Suppose $\gamma<\lambda$ and $\delta<\kappa$ and $U\in\mcU_{\gamma,\delta}$ and  $\la\la\zeta_\alpha,\eta_\alpha\ra\ra_{\alpha<\kappa}\in(\lambda\times\kappa)^\kappa$ and $\la W_\alpha\ra_{\alpha<\kappa}\in\prod_{\alpha<\kappa}\mcU_{\zeta_\alpha,\eta_\alpha}$ and $U\subseteq\bigcap_{\alpha<\kappa}W_\alpha$.  Then it suffices to show that $W_\alpha=W_\beta$ for some $\alpha<\beta<\kappa$.  Choose $V\in\mcU_{\gamma,\delta}$ such that $\closure{V}\subseteq U$.  For each $\alpha<\kappa$, choose $V_\alpha\in\mcV_{\zeta_\alpha,\eta_\alpha}$ such that $\closure{W}_\alpha\subseteq V_\alpha$; set $H_\alpha=E_{\zeta_\alpha}(W_\alpha,V_\alpha)$.  Then $E_\gamma(V,U)\subseteq\bigcap_{\alpha<\kappa}H_\alpha$.  By (\ref{enumDdense3}) and (\ref{enumDkappaop3}), $\mcD_\lambda$ is \kappaop; hence, there exists $J\in[\kappa]^{\omega_1}$ such that $H_\alpha=H_\beta$ for all $\alpha,\beta\in J$; hence, $W_\alpha\subseteq V_\beta$ for all $\alpha,\beta\in J$. If $\alpha,\beta\in J$ and $\zeta_\alpha<\zeta_\beta$, then $\bigcup\Sigma_{\zeta_\beta}\ni W_\alpha\subseteq V_\beta$, in contradiction with $V_\beta\in\mcV_{\zeta_\beta}$.  Hence, $\zeta_\alpha=\zeta_\beta$ for all $\alpha,\beta\in J$.  If $\alpha,\beta\in J$ and $\eta_\alpha<\eta_\beta$, then $\bigcup\Gamma_{\zeta_\beta,\eta_\beta}\ni W_\alpha\subseteq V_\beta$, in contradiction with $V_\beta\in\mcV_{\zeta_\beta,\eta_\beta}$.  Hence, $\eta_\alpha=\eta_\beta$ for all $\alpha,\beta\in J$.  Hence, $\{W_\alpha:\alpha\in J\}\subseteq N_{\zeta_{\min J},\eta_{\min J}}$; hence, $W_\alpha=W_\beta$ for some $\alpha<\beta<\kappa$.
\end{proof}

\begin{lemma}\label{LEMgendyadwpiw}
Let $\kappa$ be an uncountable regular cardinal; let $X$ be a compactum such that $\weight{X}\geq\kappa$ and $X$ is a continuous image of a product of compacta each with weight less than $\kappa$.  Then $\piweight{X}=\weight{X}$.
\end{lemma}
\begin{proof}
It suffices to prove that $\piweight{X}\geq\kappa$.  Seeking a contradiction, suppose $\mcA$ is a $\pi$\nbd-base of $X$ of size less than $\kappa$.  Let $\la X_i\ra_{i\in I}$ be a sequence of compacta each with weight less than $\kappa$ and let $h$ be a continuous surjection from $\prod_{i\in I}X_i$ to $X$.  Choose $M\elemsub H_\theta$ such that $\mcA\cup\{C(X),h,\la C(X_i)\ra_{i\in I}\}\subseteq M$ and $\card{M}=\card{\mcA}$.  Choose $p\in M\cap\prod_{i\in I}X_i$ and set $Y=\{q\in\prod_{i\in I}X_i:p\restrict(I\setminus M)=q\restrict(I\setminus M)\}$.  Then it suffices to show that $h[Y]=X$, for that implies $\kappa\leq\weight{X}\leq\weight{Y}<\kappa$.  Seeking a contradiction, suppose $h[Y]\not=X$.  Then there exists $U\in\mcA$ such that $U\cap h[Y]=\emptyset$.  By elementarity, there exists $\sigma\in[I\cap M]^{<\omega}$ and $\la V_i\ra_{i\in\sigma}$ such that $V_i$ is a nonempty open subset of $X_i$ for all $i\in\sigma$, and $\bigcap_{i\in\sigma}\inv{\pi_i}V_i\subseteq\inv{h}U$.  Hence, $Y\cap\bigcap_{i\in\sigma}\inv{\pi_i}V_i\not=\emptyset$, in contradiction with $U\cap h[Y]=\emptyset$.
\end{proof}

\begin{definition}
Given any cardinal $\kappa$, set $\log\kappa=\min\{\lambda:2^\lambda\geq\kappa\}$.
\end{definition}

\begin{lemma}\label{LEMgendyadwpichar}
Let $\kappa$ be an uncountable regular cardinal; let $X$ be a compactum such that $\weight{X}\geq\kappa$ and $X$ is a continuous image of a product of compacta each with weight less than $\kappa$.  Then $\pi\character{X}=\weight{X}$.
\end{lemma}
\begin{proof}
Let $\la X_i\ra_{i\in I}$ be a sequence of compacta each with weight less than $\kappa$ and let $h$ be a continuous surjection from $\prod_{i\in I}X_i$ to $X$.  For any space $Y$, we have $\piweight{Y}=\pi\character{Y}\density{Y}$.  Hence, $\weight{X}=\piweight{X}=\pi\character{X}\density{X}$ by Lemma~\ref{LEMgendyadwpiw}; hence, \wma\ $\density{X}=\weight{X}$.  Arguing as in the proof of Lemma~\ref{LEMgendyadwpiw}, if $\mcA$ is a $\pi$\nbd-base of $X$ and $\mcA\cup\{C(X),h,\la C(X_i)\ra_{i\in I}\}\subseteq M\elemsub H_\theta$, then $X$ is a continuous image of $\prod_{i\in I\cap M}X_i$; hence, \wma\ $\card{I}=\piweight{X}$.  By 5.5 of \cite{juhasz}, $\density{X}\leq\density{\prod_{i\in I}X_i}\leq\kappa\cdot\log\card{I}$.  By 2.37 of \cite{juhasz}, $\density{Y}\leq{\pi\character{Y}}^{\cell{Y}}$ for all $T_3$ non\nbd-discrete spaces $Y$.  Since $\kappa$ is a caliber of $X_i$ for all $i\in I$, it is also a caliber of $X$; hence, $\card{I}=\piweight{X}=\density{X}\leq{\pi\character{X}}^\kappa$; hence, $\log\card{I}\leq\kappa\cdot\pi\character{X}$.  Therefore, $\weight{X}=\density{X}\leq\kappa\cdot\pi\character{X}$; hence, \wma\ $\weight{X}=\kappa$.

Let $\la U_\alpha\ra_{\alpha<\kappa}$ enumerate a base of $X$.  For each $\alpha<\kappa$, choose $p_\alpha\in U_\alpha$.  Since $\density{X}=\weight{X}=\kappa$, there is no $\alpha<\kappa$ such that $\{p_\beta:\beta<\alpha\}$ is dense in $X$.  Since $\kappa$ is a caliber of $X$, we may choose $p\in X\setminus\bigcup_{\alpha<\kappa}\closure{\{p_\beta:\beta<\alpha\}}$.  It suffices to show that $\pi\character{p,X}=\kappa$.  Seeking a contradiction, suppose $\pi\character{p,X}<\kappa$.  Then there exists $\alpha<\kappa$ such that $\{U_\beta:\beta<\alpha\}$ contains a local $\pi$\nbd-base at $p$; hence, $p\in\closure{\{p_\beta:\beta<\alpha\}}$, in contradiction with how we chose $p$.
\end{proof}

\begin{theorem}\label{THMhomoggendyadic}
Let $\la X_i\ra_{i\in I}$ be a sequence of compacta; let $X$ be a homogeneous compactum; let $h\colon\prod_{i\in I}X_i\rightarrow X$ be a continuous surjection.  If there is a regular cardinal $\kappa$ such that $\weight{X_i}<\kappa\leq\weight{X}$ for all $i\in I$, then every base of $X$ contains a \cardop{(\sup_{i\in I}\weight{X_i})}\ base.  Otherwise, $\weight{X}\leq\sup_{i\in I}\weight{X_i}$ and every base of $X$ contains a \cardop{(\weight{X}^+)}\ base.
\end{theorem}
\begin{proof}
The latter case is a trivial application of Proposition~\ref{PROwnet}.  In the former case, Lemma~\ref{LEMgendyadwpichar} implies $\pi\character{p,X}=\weight{X}$ for all $p\in X$; apply Theorem~\ref{THMgendyadicweightpichar}.
\end{proof}

Every known homogeneous compactum is a continuous image of a product of compacta each with weight at most $\cardc$; hence, Theorem~\ref{THMhomoggendyadic} provides a uniform justification for our observation that all known homogeneous compacta have Noetherian type at most $\cardc^+$.  Analogously, since every known homogeneous compactum is such a continuous image, it has $\cardc^+$ among its calibers; hence, it has cellularity at most $\cardc$.

Let us now turn to the spectrum of Noetherian types of dyadic compacta and a proof of Theorem~\ref{THMdyadicspectrum}.

\begin{theorem}\label{THM2kappa2mu}
Let $\kappa$ and $\lambda$ be infinite cardinals such that $\lambda<\kappa$.  Let $X$ be the discrete sum of $2^\kappa$ and $2^\lambda$.  Let \/\/$Y$\! be the quotient space induced by collapsing $\la 0\ra_{\alpha<\kappa}$ and $\la 0\ra_{\alpha<\lambda}$ to a single point $p$.  If $\lambda<\cf\kappa$, then $\ow{Y}=\kappa^+$.  If $\lambda\geq\cf\kappa$, then $\ow{Y}=\kappa$.
\end{theorem}
\begin{proof}
Clearly $\character{p,Y}=\kappa$ and $\pi\character{p,Y}=\lambda$.  Hence, if $\lambda<\cf\kappa$, then $\kappa^+\leq\ow{Y}\leq\weight{Y}^+=\kappa^+$ by Proposition~\ref{PROpicharcfchar}.  Suppose $\lambda\geq\cf\kappa$.  We still have $\kappa\leq\ow{Y}$ by Proposition~\ref{PROpicharcfchar}, so it suffices to construct a \kappaop\ base of $Y$.  Let $\sim$ be the equivalence relation such that $Y=X/\!\sim$.  In building a base of $Y$, we proceed in the canonical way when away from $p$: for each $\mu\in\{\kappa,\lambda\}$, set
\begin{equation*}
\mcA_\mu=\{\{x\in 2^\mu:\eta\subseteq x\}/\!\sim\ :\eta\in\Fn{\mu}{2}\text{ and }\inv{\eta}\{1\}\not=\emptyset\}.
\end{equation*}
Choose $f_0\colon\kappa\to\cf\kappa$ such that for all $\alpha<\cf\kappa$ the preimage $\inv{f_0}\{\alpha\}$ is bounded in $\kappa$.  Define $f\colon[\kappa]^{<\omega}\to\cf\kappa$ by $f(\sigma)=f_0(\sup\sigma)$ for all $\sigma\in[\kappa]^{<\omega}$.  Choose $g_0\colon\lambda\to\cf\kappa$ such that for all $\alpha<\cf\kappa$ the preimage $\inv{g_0}\{\alpha\}$ is unbounded in $\lambda$.  Define $g\colon[\lambda]^{<\omega}\to\cf\kappa$ by $g(\sigma)=g_0(\sup\sigma)$ for all $\sigma\in[\lambda]^{<\omega}$.  Set
\begin{multline*}
\mcA_p=\bigcup_{\alpha<\cf\kappa}\Bigl\{\bigl(\{x\in 2^\kappa:x[\sigma]=\{0\}\}\cup\{x\in 2^\lambda:x[\tau]=\{0\}\}\bigr)/\!\sim\ :\\
\la\sigma,\tau\ra\in\inv{f}\{\alpha\}\times\inv{g}\{\alpha\}\Bigr\}.
\end{multline*}

Set $\mcA=\mcA_\kappa\cup\mcA_\lambda\cup\mcA_p$.  Let us show that $\mcA$ is a \kappaop\ base of $Y$.  The only nontrivial aspect of showing that $\mcA$ is a base of $Y$ is verifying that $\mcA_p$ is a local base at $p$.  Suppose $U$ is an open neighborhood of $p$.  Then there exist $\sigma\in[\kappa]^{<\omega}$ and $\tau\in[\lambda]^{<\omega}$ such that 
\begin{equation*}
\bigl(\{x\in 2^\kappa:x[\sigma]=\{0\}\}\cup\{x\in 2^\lambda:x[\tau]=\{0\}\}\bigr)/\!\sim\subseteq U.
\end{equation*}
Choose $\alpha<\lambda$ such that $\sup\tau<\alpha$ and $g_0(\alpha)=f(\sigma)$.  Set $\tau'=\tau\cup\{\alpha\}$ and
\begin{equation*}
V=\bigl(\{x\in 2^\kappa:x[\sigma]=\{0\}\}\cup\{x\in 2^\lambda:x[\tau']=\{0\}\}\bigr)/\!\sim.
\end{equation*}
Then $V\subseteq U$ and $V\in\mcA_p$ because $f(\sigma)=g(\tau')$.  Thus, $\mcA$ is a base of $Y$.

Let us show that $\mcA$ is \kappaop.  Suppose $U,V\in\mcA$ and $U\subseteq V$.  If $U\in\mcA_\kappa$, then, fixing $U$, there are only finitely possibilities for $V$ in $\mcA_\kappa$; the same is true if $\kappa$ is replaced by $\lambda$ or $p$.  Hence, we may assume $U\in\mcA_i$ and $V\in\mcA_j$ for some $\{i,j\}\in[\{\kappa,\lambda,p\}]^2$.  Since no element of $\mcA_p$ is a subset of an element of $\mcA_\kappa\cup\mcA_\lambda$, we have $i\not=p$.  Hence, there exists $\eta\in\Fn{i}{2}$ such that $U=\{x\in 2^i:\eta\subseteq x\}/\!\sim$.  Since $\bigcup\mcA_\kappa\cap\bigcup\mcA_\lambda=\emptyset$, we have $j=p$.  Hence, there exist $\sigma\in[\kappa]^{<\omega}$ and $\tau\in[\lambda]^{<\omega}$ such that 
\begin{equation*}
V=\bigl(\{x\in 2^\kappa:x[\sigma]=\{0\}\}\cup\{x\in 2^\lambda:x[\tau]=\{0\}\}\bigr)/\!\sim.
\end{equation*}
If $i=\kappa$, then $\sigma\subseteq\inv{\eta}\{0\}$; hence, fixing $U$, there are only finitely many possibilities for $\sigma$, and at most $\lambda$\nbd-many possibilities for $\tau$.  If $i=\lambda$, then $\tau\subseteq\inv{\eta}\{0\}$; hence, fixing $U$, there are only finitely many possibilities for $\tau$, and at most $\card{\sup\inv{f_0}\{g(\tau)\}}^{<\omega}$\nbd-many possibilities for $\sigma$ given $\tau$.  Thus, there are fewer than $\kappa$\nbd-many possibilities for $V$ given $U$.  Thus, $\mcA$ is \kappaop.
\end{proof}

\begin{corollary}\label{COR2kappa2muspectrum}
If $\kappa$ is a cardinal of uncountable cofinality, then there is a totally disconnected dyadic compactum with Noetherian type $\kappa^+$.  If $\kappa$ is a singular cardinal, then there is a totally disconnected dyadic compactum with Noetherian type $\kappa$.
\end{corollary}
\begin{proof}
For the first case, apply Theorem~\ref{THM2kappa2mu} with $\lambda=\omega$.  For the second case, apply Theorem~\ref{THM2kappa2mu} with $\lambda=\cf\kappa$.
\end{proof}

Combining the above corollary with the following theorem (and a trivial example like $\ow{2^\omega}=\omega$) immediately proves Theorem~\ref{THMdyadicspectrum}.

\begin{theorem}
Let $X$ be a dyadic compactum with base $\mcA$ consisting only of cozero sets.  If $\ow{X}\leq\omega_1$, then $\mcA$ contains an \omegaop\ base of $X$.  Hence, no dyadic compactum has Noetherian type $\omega_1$.
\end{theorem}
\begin{proof}
Let $\mcQ$ be an \cardop{\omega_1}\ base of $X$ of size $\weight{X}$.  Import all the notation from the proof of Lemma~\ref{LEMweightpicharcozero} verbatim, except that we require $\la M_\alpha\ra_{\alpha<\kappa}$ to be an $\omega_1$\nbd-approximation sequence in $\la H_\theta,\in,\mcF,h,\mcQ\ra$.  Then $\mcU$ is an \omegaop\ subset of $\mcA$ as before.  On the other hand, $\mcV_\alpha/\mcF_\alpha$ is not necessarily a base of $X/\mcF_\alpha$ for all $\alpha<\kappa$. However, we will show that $\mcU$ is still a base of $X$.  In doing so, we will repeatedly use the fact that if $U,\mcQ\in M\elemsub H_\theta$ and $U$ is a nonempty open subset of $X$, then all supersets of $U$ in $\mcQ$ are in $M$ because $\{V\in\mcQ:U\subseteq V\}$ is a countable element of $M$.

Suppose $q\in Q\in\mcQ$.  Then it suffices to find $U\in\mcU$ such that $q\in U\subseteq Q$.  Let $\beta$ be the least $\alpha<\kappa$ such that there exists $A\in\mcA_\alpha$ satisfying $q\in A\subseteq\closure{A}\subseteq Q$.  Fix such an $A\in\mcA_\beta$.  
For each $p\in\closure{A}$, choose $\la A_p, Q_p\ra\in\mcA\times\mcQ$ such that $p\in A_p\subseteq  Q_p\subseteq\closure{Q}_p\subseteq Q$.  Since $M_\beta\ni A\subseteq Q\in\mcQ$, we have $Q\in M_\beta$.  Hence, by elementarity, \wma\ there exists $\sigma\in\bigl[\,\closure{A}\,\bigr]^{<\omega}$ such that $\la\la A_p,Q_p\ra\ra_{p\in\sigma}\in M_\beta$ and $\closure{A}\subseteq\bigcup_{p\in\sigma}A_p$.  Choose $p\in\sigma$ such that $q\in A_p$.  Suppose $Q_p\not\in\bigcup\Sigma_\beta$.  Then all nonempty open subsets of $Q_p$ are also not in $\bigcup\Sigma_\beta$; hence, there exist $U\in\mcU_\beta$ and $V\in\mcV_\beta$ such that $q/\mcF_\beta\subseteq U\subseteq V\subseteq A_p\subseteq Q$.  Therefore, \wma\ $Q_p\in\bigcup\Sigma_\beta$.

Choose $\alpha<\beta$ such that $Q_p\in M_\alpha$.  Then $Q\in M_\alpha$ because $Q_p\subseteq Q$.  Hence, there exists $\tau\in[\mcA_\alpha]^{<\omega}$ such that $\closure{Q}_p\subseteq\bigcup\tau\subseteq\closure{\bigcup\tau}\subseteq Q$.  Choose $W\in\tau$ such that $q\in W$.  Then $q\in W\subseteq\closure{W}\subseteq Q$, in contradiction with the minimality of $\beta$.  Thus, $\mcU$ is a base of $X$.
\end{proof}

\begin{question}
If $\kappa$ is an singular cardinal with cofinality $\omega$, then is there a dyadic compactum with Noetherian type $\kappa^+$?  Is there a dyadic compactum with weakly inaccessible Noetherian type?
\end{question}

We note that the spectrum of Noetherian types of all compacta is trivial.

\begin{theorem}\label{THMowkappabool}
Let $\kappa$ be a regular uncountable cardinal.  Then there exists a totally disconnected compactum $X$ such that $\ow{X}=\kappa$ and $X$ has a $P_\kappa$\nbd-point.
\end{theorem}
\begin{proof}
Let $X$ be the closed subspace of $2^\kappa$ consisting of all $f\in 2^\kappa$ for which $f(\alpha)=0$ or $f[\alpha]=\{1\}$ for all odd $\alpha<\kappa$.  First, let us show that $X$ has a \kappaop\ base.  For each $\sigma\in\Fn{\kappa}{2}$, set $U_\sigma=\{f\in X:f\supseteq\sigma\}$.  Let $E$ denote the set of $\sigma\in\Fn{\kappa}{2}$ for which $\sup\dom\sigma$ is even and $U_\sigma\not=\emptyset$.  Set $\mcA=\{U_\sigma:\sigma\in E\}$, which is clearly a base of $X$.  Let us show that $\mcA$ is \kappaop.  Suppose $\sigma,\tau\in E$ and $U_\sigma\subseteq U_\tau$.  If $\sup\dom\sigma<\sup\dom\tau$, then for each $f\in U_\sigma$ the sequence
\begin{equation*}
(f\restrict\sup\dom\tau)\cup\{\la\sup\dom\tau,1-\tau(\sup\dom\tau)\ra\}\cup\{\la\beta,0\ra:\sup\dom\tau<\beta<\kappa\}
\end{equation*}
is in $U_\sigma\setminus U_\tau$, which is absurd.  Hence, $\sup\dom\tau\leq\sup\dom\sigma$; hence, there are fewer than $\kappa$\nbd-many possibilities for $\tau$ given $\sigma$.  Thus, $\mcA$ is \kappaop.

Finally, it suffices to show that $\la 1\ra_{\alpha<\kappa}$ is a $P_\kappa$\nbd-point of $X$, for a $P_\kappa$\nbd-point must have local Noetherian type at least $\kappa$.  For each $\alpha<\kappa$, set $\sigma_\alpha=\{\la 2\alpha+1,1\ra\}$.  Then $\{U_{\sigma_\alpha}:\alpha<\kappa\}$ is a local base at $\la 1\ra_{\alpha<\kappa}$.  Moreover, $U_{\sigma_\alpha}\supsetneq U_{\sigma_\beta}$ for all $\alpha<\beta<\kappa$.  Since $\kappa$ is regular, it follows that $\la 1\ra_{\alpha<\kappa}$ is a $P_\kappa$\nbd-point.
\end{proof}

\begin{corollary}
Every infinite cardinal is the Noetherian type of some totally disconnected compactum.
\end{corollary}
\begin{proof}
By Lemma~\ref{LEMcpctmetricbase}, all totally disconnected metric compacta have Noetherian type $\omega$.  By Theorem~\ref{THMowkappabool}, if $\kappa$ is a regular uncountable cardinal, then there is a totally disconnected compactum $X$ with Noetherian type $\kappa$.  If $\kappa$ is a singular cardinal, then there is a totally disconnected dyadic compactum with Noetherian type $\kappa$ by Corollary~\ref{COR2kappa2muspectrum}.
\end{proof}

\section{Reflecting cones}\label{SECreflectcone}

Notice that the proofs of Theorems~\ref{THMalmostomegaopsubsetsdyadic} and \ref{THMweightpichar} only indirectly use the hypothesis that $X$ is dyadic: if $X$ is merely a continuous Hausdorff image of the Stone space of a boolean algebra satisfying the conclusion of Lemma~\ref{LEMfreeboolreflect}, then it is routine to check that the proofs of Theorems~\ref{THMalmostomegaopsubsetsdyadic} and \ref{THMweightpichar} are still valid.  This prompts the question of whether there are nondyadic compacta $X$ for which these proofs apply.  Equivalently, is there a boolean algebra $B$ such that $B$ is not a subalgebra of a free boolean algebra but the conclusion of Lemma~\ref{LEMfreeboolreflect} holds for $B$?  We show that the answer is no if $\card{B}\leq\omega_1$, and present some partial results for larger $B$.

\begin{definition}
Let $B$ be a boolean algebra.  We say that $B$ \emph{reflects cones} to $A$ if $A$ is a subalgebra of $B$ and, for all $q\in B$, there exists $r\in A$ such that for all $p\in A$ we have $p\geq q$ if and only if $p\geq r$.  Let $\la H_\theta,\ldots\ra$ denote an expansion of the $\{\in\}$\nbd-structure $H_\theta$ to some $\mcL$\nbd-structure for some countable language $\mcL$.  We say that $B$ reflects cones in $\la H_\theta,\ldots\ra$ if $B$ reflects cones to $B\cap M$ for all $M\elemsub\la H_\theta,\ldots\ra$.  If $n<\omega$, then we say that $B$ $n$\nbd-reflects cones in $\la H_\theta,\ldots\ra$ if, for all $\Sigma$ of size at most $n$ satisfying (\ref{enumsigmaelemsub}), (\ref{enumepsilonchain}), and (\ref{enumdesccard}) of Lemma~\ref{LEMapprox}, $B$ reflects cones to the subalgebra generated by $B\cap\bigcup\Sigma$.  We say that $B$ $\omega$\nbd-reflects cones in $\la H_\theta,\ldots\ra$ if $B$ $n$\nbd-reflects cones in $\la H_\theta,\ldots\ra$ for all $n<\omega$.  We say that $B$ reflects cones ($n$\nbd-reflects cones, $\omega$\nbd-reflects cones) if $B$ reflects cones ($n$\nbd-reflects cones, $\omega$\nbd-reflects cones) in some $\la H_\theta,\ldots\ra$ for all sufficiently large $\theta$.
\end{definition}

Note that $1$\nbd-reflecting cones is equivalent to reflecting cones.

\begin{lemma}\label{LEMreflectonemore}
Suppose $B$ is boolean algebra with element $c$ and subalgebra $A$ such that $B$ reflects cones to $A$.  Let $C$ denote the subalgebra of $B$ generated by $A\cup\{c\}$.  Then $B$ reflects cones to $C$.
\end{lemma}
\begin{proof}
Let $q\in B$. By hypothesis, there exist $r_0,r_1\in A$ such that for all $p\in A$ we have $p\geq q\land c$ if and only if $p\geq r_0$ and $p\geq q\land c'$ if and only if $p\geq r_1$.  In particular, $r_0\geq q\land c$ and $r_1\geq q\land c'$.  Set $r=(r_0\lor c')\land(r_1\lor c)\in C$.  Then $r\geq((q\land c)\lor c')\land((q\land c')\lor c)=q$.  Suppose $p\in C$ and $p\geq q$.  Then it suffices to show that $p\geq r$.  There exist $p_0,p_1\in A$ such that $p=(p_0\lor c)\land(p_1\lor c')$.  Then $p_0\lor c\geq q$; hence, $p_0\geq q\land c'$; hence, $p_0\geq r_1$; hence, $p_0\lor c\geq r_1\lor c$.  By symmetry, $p_1\lor c'\geq r_0\lor c'$.  Hence, $p\geq r$.
\end{proof}

\begin{lemma}\label{LEMembedonemore}
Suppose $B$ is boolean algebra with element $c$ and subalgebra $A$ such that $B$ reflects cones to $A$.  Further suppose $f$ is an embedding of $A$ into a free boolean algebra.  Let $C$ denote the subalgebra of $B$ generated by $A\cup\{c\}$.  Then $f$ extends to an embedding of $C$ into a free boolean algebra.
\end{lemma}
\begin{proof}
\Wma\ $f$ embeds $A$ into a free boolean algebra $D$ such that there exists nontrivial $d\in D$ independent from the subalgebra $f[A]$.  By hypothesis, there exist $a_0,a_1\in A$ such that for all $p\in A$ we have $p\geq c$ if and only if $p\geq a_0$ and $p\geq c'$ if and only if $p\geq a_1$.  Hence, for all $p\in A$, we have $p\leq c$ if and only if $p\leq a_1'$.  In particular, $a_0\geq c\geq a_1'$.  Let $g$ be the unique homomorphism from $C$ to $D$ extending $f$ such that $g(c)=(f(a_0)\land d)\lor f(a_1)'$.  Suppose $x,y\in C$.  Then it suffices to show that $x\geq y$ if and only if $g(x)\geq g(y)$.  There exist $x_0,x_1,y_0,y_1\in A$ such that $x=(x_0\land c)\lor(x_1\land c')$ and $y=(y_0\land c)\lor(y_1\land c')$.  Then $x\geq y$ if and only if $x_0\land c\geq y_0\land c$ and $x_1\land c'\geq y_1\land c'$; likewise, $g(x)\geq g(y)$ if and only if $f(x_0)\land g(c)\geq f(y_0)\land g(c)$ and $f(x_1)\land g(c)'\geq f(y_1)\land g(c)'$.  By symmetry, it suffices to show that $x_0\land c\geq y_0\land c$ if and only if $f(x_0)\land g(c)\geq f(y_0)\land g(c)$.  Clearly $x_0\land c\geq y_0\land c$ if and only if $x_0\lor y_0'\geq c$; likewise, $f(x_0)\land g(c)\geq f(y_0)\land g(c)$ if and only if $f(x_0)\lor f(y_0)'\geq g(c)$.  By definition, $f(x_0)\lor f(y_0)'\geq g(c)$ if and only if $f(x_0)\lor f(y_0)'\geq(f(a_0)\land d)\lor f(a_1)'$, which is equivalent to $f(x_0)\lor f(y_0)'\geq f(a_0)$ because $f(a_0)\geq f(a_1)'$ and  $d$ is independent from $f[A]$.  Hence, $x_0\land c\geq y_0\land c$ if and only if $x_0\lor y_0'\geq c$ if and only if $x_0\lor y_0'\geq a_0$ if and only if $f(x_0)\lor f(y_0)'\geq f(a_0)$ if and only if $f(x_0)\lor f(y_0)'\geq g(c)$ if and only if $f(x_0)\land g(c)\geq f(y_0)\land g(c)$.
\end{proof}

\begin{theorem}\label{THMomegareflectconesembed}
Suppose $B$ is a boolean algebra that $\omega$\nbd-reflects cones.  Then $B$ is a subalgebra of a free boolean algebra.
\end{theorem}
\begin{proof}
Let $\la M_\alpha\ra_{\alpha<\card{B}}$ be an $\omega_1$\nbd-approx\-imation sequence in $\la H_\theta,\in,B,\ldots\ra$ where $B$ $\omega$\nbd-reflects cones in $\la H_\theta,\in,B,\ldots\ra$.  Let $\la\Sigma_\alpha\ra_{\alpha\leq\card{B}}$ be as in Lemma~\ref{LEMapprox}.  For each $\alpha\leq\card{B}$, let $B_\alpha$ be the subalgebra of $B$ generated by $B\cap\bigcup\Sigma_\alpha$.  Trivially, we may choose an embedding $f_0$ of $B_0$ into a free boolean algebra.  Suppose $\alpha<\card{B}$ and we have an embedding $f_\alpha$ of $B_\alpha$ into a free boolean algebra.  Let $B\cap M_\alpha=\{b_n:n<\omega\}$.  For each $n<\omega$, let $A_n$ be the subalgebra of $B$ generated by $B_\alpha\cup\{b_m:m<n\}$; by repeated application of Lemma~\ref{LEMreflectonemore}, $B$ reflects cones to each $A_n$.  Set $g_0=f_\alpha$.  Suppose $n<\omega$ and $g_n$ is an embedding of $A_n$ into a free boolean algebra.  By Lemma~\ref{LEMembedonemore}, there is an extension $g_{n+1}$ of $g_n$ embedding $A_{n+1}$ into a free boolean algebra.  Set $f_{\alpha+1}=\bigcup_{n<\omega}g_n$.  Then $f_{\alpha+1}$ is an extension of $f_\alpha$ embedding $B_{\alpha+1}$ into a free boolean algebra.  For limit ordinals $\alpha\leq\card{B}$, set $f_\alpha=\bigcup_{\beta<\alpha}f_\beta$.  Then $f_{\card{B}}$ embeds $B$ into a free boolean algebra.
\end{proof}

\begin{corollary}
Let $n<\omega$ and $B$ be a boolean algebra of size at most $\omega_n$.  If $B$ $n$\nbd-reflects cones, then $B$ is a subalgebra of a free boolean algebra.
\end{corollary}
\begin{proof}
If $\Sigma$ satisfies (\ref{enumsigmaelemsub}), (\ref{enumepsilonchain}), and (\ref{enumdesccard}) of Lemma~\ref{LEMapprox} and $\card{\Sigma}>n$, then there exists $N\in\Sigma$ such that $\card{B}\subseteq N$.  \Wma\ $B\in N$; hence, $B\subseteq N$; hence, $B\cap\bigcup\Sigma=B$.  Thus, $B$ $\omega$\nbd-reflects cones.
\end{proof}

\begin{question}
Is there a boolean algebra that reflects cones but is not a subalgebra of a free boolean algebra?
\end{question}

\section{More on local Noetherian type}\label{SEChomogloc}

In this section, we find two sufficient conditions for a compactum to have a point with an \omegaop\ local base.  The first of these conditions will be used to prove Theorem~\ref{THMgchocharcell}.  We also present some related results about local bases in terms of Tukey reducibility.

\begin{definition}
Given cardinals $\lambda\geq\kappa\geq\omega$ and a subset $E$ in a space $X$, a \emph{local} $\la\lambda,\kappa\ra$\nbd-\emph{splitter} at $E$ is a set $\mcU$ of $\lambda$\nbd-many open neighborhoods of $E$ such that $E$ is contained in the interior of $\bigcap\mcV$ for any $\mcV\in[\mcU]^\kappa$.  If $p\in X$, then we call a local $\la\lambda,\kappa\ra$\nbd-splitter at $\{p\}$ a local $\la\lambda,\kappa\ra$\nbd-splitter a $p$.
\end{definition}

\begin{theorem}\label{THMpicharomegalocalsplitter}
Suppose $X$ is a compactum and $\omega_1\leq\kappa=\min_{p\in X}\pi\character{p,X}$.  Then there is a local $\la\kappa,\omega\ra$\nbd-splitter at some $p\in X$.
\end{theorem}
\begin{proof}
Given any map $f$, let $\prod f$ denote $\{\la x_i\ra_{i\in\dom f}:\forall i\in\dom f\ \, x_i\in f(i)\}$.  Given any infinite open family $\mcE$, let $\Phi(\mcE)$ denote the set of $\la\sigma,\Gamma\ra\in[\mcE]^{<\omega}\times([\mcE]^\omega)^{<\omega}$ for which every $\tau\in\prod\Gamma$ satisfies $\bigcap\sigma\subseteq\closure{\bigcup\ran\tau}$.  Then $\Phi(\mcE)=\emptyset$ always implies $\mcE$ is \omegaop\ and centered.  

Let $\mcR$ denote the set of nonempty regular open subsets of $X$.  Choose $\la W_n\ra_{n<\omega}\in\mcR^\omega$ such that $\closure{W}_{n+1}\subsetneq W_n\not=X$ for all $n<\omega$.  Let $\Omega$ denote the class of transfinite sequeces $\la\la U_\alpha,V_\alpha\ra\ra_{\alpha<\eta}$ of elements of $\mcR^2$ satisfying the following.
\begin{enumerate}
\item\label{enumOmegainit} $\eta\geq\omega$ and $\la\la U_n,V_n\ra\ra_{n<\omega}=\la\la W_{n+1},W_n\ra\ra_{n<\omega}$.
\item\label{enumOmegasubsetpair} $\closure{U}_\alpha\subseteq V_\alpha$ for all $\alpha<\eta$.
\item\label{enumOmegastable} $\powset{V_\alpha}\cap\Bigl\{\bigcap\sigma\setminus\closure{\bigcup\tau}:\sigma,\tau\in\bigl[\bigcup_{\beta<\alpha}\{U_\beta,V_\beta\}\bigr]^{<\omega}\Bigr\}\subseteq\{\emptyset\}$ for all $\alpha<\eta$.
\item\label{enumOmegaPhi} $\Phi\bigl(\bigcup_{\alpha<\eta}\{U_\alpha,V_\alpha\}\bigr)=\emptyset$.
\end{enumerate}

Seeking a contradiction, suppose $\eta$ is a limit ordinal and $\la\la U_\alpha,V_\alpha\ra\ra_{\alpha<\eta}\not\in\Omega$, but $\la\la U_\beta,V_\beta\ra\ra_{\beta<\alpha}\in\Omega$ for all $\alpha<\eta$.  Then (\ref{enumOmegainit}), (\ref{enumOmegasubsetpair}), and (\ref{enumOmegastable}) hold for $\la\la U_\alpha,V_\alpha\ra\ra_{\alpha<\eta}$, so there exists $\la\sigma,\Gamma\ra\in\Phi\bigl(\bigcup_{\alpha<\eta}\{U_\alpha,V_\alpha\}\bigr)$.  We may choose $i\in\dom\Gamma$ such that $\Gamma(i)\not\subseteq\bigcup_{\beta<\alpha}\{U_\beta,V_\beta\}$ for all $\alpha<\eta$.  Set $\Lambda=\Gamma\restrict(\dom\Gamma\setminus\{i\})$.  \Wma\ $\dom\Gamma$ is minimal among its possible values; hence, there exists $\tau\in\prod\Lambda$ such that $\bigcap\sigma\not\subseteq\closure{\bigcup\ran\tau}$.  Choose $\alpha<\eta$ and $W\in\Gamma(i)$ such that $\sigma\cup\ran\tau\subseteq\bigcup_{\beta<\alpha}\{U_\beta,V_\beta\}$ and $W\in\{U_\alpha,V_\alpha\}$.  Then $\bigcap\sigma\setminus\closure{\bigcup\ran\tau}\not\subseteq W$ by (\ref{enumOmegasubsetpair}) and (\ref{enumOmegastable}).  Since $W$ is regular, $\bigcap\sigma\setminus\closure{\bigcup\ran\tau}\not\subseteq\closure{W}$; hence, $\bigcap\sigma\not\subseteq\closure{W\cup\bigcup\ran\tau}$, in contradiction with $\la\sigma,\Gamma\ra\in\Phi\bigl(\bigcup_{\alpha<\eta}\{U_\alpha,V_\alpha\}\bigr)$.  Thus, $\Omega$ is closed with respect to unions of increasing chains.  

It follows from (\ref{enumOmegastable}) that $\Omega\subseteq(\mcR^2)^{<\card{\mcR}^+}$.  Moreover, $\la\la W_{n+1},W_n\ra\ra_{n<\omega}\in\Omega$.  Hence, by Zorn's Lemma, $\Omega$ has a maximal element $\la\la U_\alpha,V_\alpha\ra\ra_{\alpha<\eta}$.  Set $\mcB=\bigcup_{\alpha<\eta}\{U_\alpha,V_\alpha\}$.  Let us show that $\eta\geq\kappa$.  Suppose not.  For each $x\in X$, choose $Y_x,Z_x\in\mcR$ such that $x\in Y_x\subseteq\closure{Y}_x\subseteq Z_x$ and $Z_x$ does not contain any nonempty open set of the form $\bigcap\sigma\setminus\closure{\bigcup\tau}$ where $\sigma,\tau\in[\mcB]^{<\omega}$.  Choose $\rho\in[X]^{<\omega}$ such that $\bigcup_{x\in\rho}Y_x=X$.  Let us show that $\Phi(\mcB\cup\{Y_x,Z_x\})=\emptyset$ for some $x\in\rho$.  Seeking a contradiction, suppose $\la\sigma_x,\Gamma_x\ra\in\Phi(\mcB\cup\{Y_x,Z_x\})$ for all $x\in\rho$.  \Wma\ $\bigcup_{x\in\rho}\bigcup\ran\Gamma_x\subseteq\mcB$.  Let $\Lambda$ be a concatenation of $\{\Gamma_x:x\in\rho\}$ and set $\tau=\mcB\cap\bigcup_{x\in\rho}\sigma_i$.  Then for all $\zeta\in\prod\Lambda$ we have 
\begin{equation*}
\bigcap\tau=\bigcap_{y\in\rho}\bigcap(\sigma_y\cap\mcB)=\bigcup_{x\in\rho}\left(Y_x\cap\bigcap_{y\in\rho}\bigcap(\sigma_y\cap\mcB)\right)\subseteq\bigcup_{x\in\rho}\bigcap\sigma_x\subseteq\closure{\bigcup\ran\zeta}.
\end{equation*}
Hence, $\la\tau,\Lambda\ra\in\Phi(\mcB)$, in contradiction with (\ref{enumOmegaPhi}).  Therefore, we may choose $x\in\rho$ such that $\Phi(\mcB\cup\{Y_x,Z_x\})=\emptyset$.  But then $\la\la U_\alpha,V_\alpha\ra\ra_{\alpha<\eta+1}\in\Omega$ if we set $U_\eta=Y_x$ and $V_\eta=Z_x$, in contradiction with the maximality of $\la\la U_\alpha,V_\alpha\ra\ra_{\alpha<\eta}$.  Thus, $\eta\geq\kappa$.

Set $\mcA=\{V_\alpha:\alpha<\eta\}$.  By (\ref{enumOmegastable}), $\card{\mcA}=\card{\eta}\geq\kappa$.  Set $K=\bigcap_{\alpha<\eta}\closure{U}_\alpha$.  Then it suffices to show that $\mcA$ is a local $\la\card{\eta},\omega\ra$\nbd-splitter at some $x\in K$.  Suppose not.  Then each $x\in K$ has an open neighborhood $W_x$ that is a subset of infinitely many elements of $\mcA$.  Hence, $\Phi(\mcB\cup\{W_x\})\not=\emptyset$ for all $x\in K$.  Choose $\rho\in[K]^{<\omega}$ such that $K\subseteq\bigcup_{x\in\rho}W_x$.   Choose an open set $W$ such that $W\cup\bigcup_{x\in\rho}W_x=X$ and $\closure{W}\cap K=\emptyset$.  By compactness, $\mcB\cup\{W\}$ is not centered; hence, $\Phi(\mcB\cup\{W\})\not=\emptyset$.  Reusing our earlier concatenation argument, we have $\Phi(\mcB)\not=\emptyset$, in contradiction with (\ref{enumOmegaPhi}).  Thus, $\mcA$ is a local $\la\card{\eta},\omega\ra$\nbd-splitter at some $x\in K$.
\end{proof}

\begin{lemma}\label{LEMocharconverse}
Suppose $E$ is a subset of a space $X$ and $E$ has no finite neighborhood base.  Then $\ochar{E,X}$ is the least $\kappa\geq\omega$ for which there is a local $\la\character{p,X},\kappa\ra$\nbd-splitter at $E$.
\end{lemma}
\begin{proof}
Set $\kappa=\ochar{E,X}$ and $\lambda=\character{E,X}$.  By Lemma~\ref{LEMkappasizeposet}, $\lambda\geq\kappa$; hence, a \kappaop\ neighborhood base of $E$ (which necessarily has size $\lambda$) is a local $\la\lambda,\kappa\ra$\nbd-splitter at $E$.  To show the converse, let $\la U_\alpha\ra_{\alpha<\lambda}$ be a sequence of open neighborhoods of $E$.  Let $\{V_\alpha:\alpha<\lambda\}$ be a neighborhood base of $E$.  For each $\alpha<\lambda$, choose $W_\alpha\in\{V_\beta:\beta<\lambda\}$ such that $W_\alpha\subseteq U_\alpha\cap V_\alpha$.  Then $\{W_\alpha:\alpha<\lambda\}$ is a neighborhood base of $E$.  Let $\mu<\kappa$.  Then there exist $\alpha<\lambda$ and $I\in[\lambda]^\mu$ such that $W_\alpha\subseteq\bigcap_{\beta\in I}W_\beta$.  Hence, $E$ is contained in the interior of $\bigcap_{\beta\in I}U_\beta$.  Hence, $\{U_\alpha:\alpha<\lambda\}$ is not a local $\la\lambda,\mu\ra$\nbd-splitter at $E$.
\end{proof}

\begin{proof}[Proof of Theorem~\ref{THMocharpicharchar}]
\Wma\ $\character{X}\geq\omega_1$.  By Theorem~\ref{THMpicharomegalocalsplitter}, there is a local $\la\character{X},\omega\ra$\nbd-splitter at some $p\in X$.  By Lemma~\ref{LEMocharconverse}, $\ochar{p,X}=\omega$.
\end{proof}

\begin{proof}[Proof of Theorem~\ref{THMgchocharcell}]
Let $X$ be a homogeneous compactum.  By a result of Arhan\-gel${}^\prime$ski\u\i\ (see 1.5 of \cite{arhangelskii}), $\card{Y}\leq 2^{\pi\character{Y}\cell{Y}}$ for all homogeneous spaces $Y$.  Since $\card{X}=2^{\character{X}}$ by Arhangel${}^\prime$ski\u\i's Theorem and the \v Cech\nbd-Pospi\v sil Theorem, we have $\character{X}\leq\pi\character{X}\cell{X}$ by GCH.  If $\pi\character{X}=\character{X}$, then $\ochar{X}=\omega$ by Theorem~\ref{THMocharpicharchar}.  Hence, \wma\ $\pi\character{X}<\character{X}$; hence, $\ochar{X}\leq\character{X}\leq\cell{X}$ by Theorem~\ref{THMcharpibounds}.
\end{proof}

\begin{example}
Consider $2^{\omega_1}$ ordered lexicographically.  Every point in this space has character and local Noetherian type $\omega_1$, and some but not all points have $\pi$\nbd-character $\omega$.
\end{example}

\begin{definition}[Tukey~\cite{tukey}]
Given two quasiorders $P$ and $Q$, we say $f$ is a \emph{Tukey} map from $P$ to $Q$ and write $f\colon P\leq_T Q$ if $f$ is a map from $P$ to $Q$ such that all preimages of bounded subsets of $Q$ are bounded in $P$.  We say that $P$ is \emph{Tukey reducible} to $Q$ and write $P\leq_T Q$ if there exists $f\colon P\leq_T Q$.  We say that $P$ and $Q$ are \emph{Tukey equivalent} and write $P\equiv_T Q$ if $P\leq_T Q\leq_T P$.
\end{definition}

Tukey showed that two directed sets are Tukey equivalent if and only if they embed as cofinal subsets of a common directed set.  In particular, any two local bases at a common point in a topological space are Tukey equivalent.  Another, easily checked fact is thats $P\leq_T[\cf P]^{<\omega}$ for every directed set $P$.  Also, $[\kappa]^{<\omega}\leq_T[\lambda]^{<\omega}$ if $\kappa\leq\lambda$.

\begin{lemma}\label{LEMlocsplittukey}
Suppose $\kappa\geq\omega$ and $E$ is a subset of a space $X$ with a local $\la\kappa,\omega\ra$\nbd-splitter at $E$.  Then $\la[\kappa]^{<\omega},\subseteq\ra\leq_T\la\mcA,\supseteq\ra$ for every neighborhood base $\mcA$ of $E$.
\end{lemma}
\begin{proof}
Let $\mcU$ be a local $\la\kappa,\omega\ra$\nbd-splitter at $E$.  Let $\mcN$ be the set of open neighborhoods of $E$.  Then $\mcN$ is Tukey equivalent to every neighborhood base of $E$ (with respect to $\supseteq$), so it suffices to show that $[\mcU]^{<\omega}\leq_T\la\mcN,\supseteq\ra$.  Define $f\colon[\mcU]^{<\omega}\rightarrow\mcN$ by $f(\sigma)=\bigcap\sigma$ for all $\sigma\in[\mcU]^{<\omega}$.  Then, for all $N\in\mcN$, we have $\card{\inv{f}\upset N}<\omega$ because $\mcU$ is a local $\la\kappa,\omega\ra$\nbd-splitter; whence, $\inv{f}\upset N$ is bounded in $[\mcU]^{<\omega}$.  Thus, $f\colon[\mcU]^{<\omega}\leq_T\la\mcN,\supseteq\ra$.
\end{proof}

\begin{theorem}\label{THMpichartukey}
Suppose $X$ is a compactum and $\omega_1\leq\kappa=\min_{p\in X}\pi\character{p,X}$.  Then, for some $p\in X$, every local base $\mcA$ at $p$ satisfies $\la[\kappa]^{<\omega},\subseteq\ra\leq_T\la\mcA,\supseteq\ra$.
\end{theorem}
\begin{proof}
Combine Theorem~\ref{THMpicharomegalocalsplitter} and Lemma~\ref{LEMlocsplittukey}.
\end{proof}

\begin{lemma}\label{LEMtukeylocntlocsplit}
Suppose $E$ is a subset of a space $X$ and $E$ has no finite neighborhood base.  Then the following are equivalent.
\begin{enumerate}
\item\label{enumeqntloc} $\ochar{E,X}=\omega$.
\item\label{enumeqlocsplit} There is a local $\la\character{E,X},\omega\ra$\nbd-splitter at $E$.
\item\label{enumeqtukey} Every neighborhood base $\mcA$ of $E$ satisfies $\la[\character{E,X}]^{<\omega},\subseteq\ra\equiv_T\la\mcA,\supseteq\ra$.
\end{enumerate}
\end{lemma}
\begin{proof}
By Lemma~\ref{LEMocharconverse}, (\ref{enumeqntloc}) and (\ref{enumeqlocsplit}) are equivalent.  Let $\mcB$ be a neighborhood base of $E$ of size $\character{E,X}$.  By Lemma~\ref{LEMlocsplittukey}, (\ref{enumeqlocsplit}) implies $[\character{E,X}]^{<\omega}\leq_T\la\mcA,\supseteq\ra\equiv_T\la\mcB,\supseteq\ra\leq_T[\character{E,X}]^{<\omega}$ for every neighborhood base $\mcA$ of $E$.  Thus, (\ref{enumeqlocsplit}) implies (\ref{enumeqtukey}).  Finally, suppose $\mcA$ is a neighborhood base of $E$ and $[\character{E,X}]^{<\omega}\equiv_T\la\mcA,\supseteq\ra$.  Then $[\character{E,X}]^{<\omega}$ and $\la\mcA,\supseteq\ra$ embed as cofinal subsets of a common directed set.  Hence, $\la\mcA,\subseteq\ra$ is almost \omegaop\ by Lemma~\ref{LEMmutuallydense}.  Hence, $\mcA$ contains an \omegaop\ neighborhood base of $E$.  Thus, (\ref{enumeqtukey}) implies (\ref{enumeqntloc}).
\end{proof}

\begin{theorem}
Suppose $X$ is an infinite homogeneous compactum and $\pi\character{X}=\character{X}$.  Then, for all $p\in X$ and for all local bases $\mcA$ at $p$, we have $\la\mcA,\supseteq\ra\equiv_T\la[\character{X}]^{<\omega},\subseteq\ra$.
\end{theorem}
\begin{proof}
Combine Theorem~\ref{THMocharpicharchar} and Lemma~\ref{LEMtukeylocntlocsplit}.
\end{proof}

\begin{definition}
Given $n<\omega$ and ordinals $\alpha,\beta_0,\ldots,\beta_n$, let $\alpha\rightarrow(\beta_0,\ldots,\beta_n)$ denote the proposition that for all $f\colon[\alpha]^2\rightarrow n+1$ there exist $i\leq n$ and $H\subseteq\alpha$ such that $f[[H]^2]=\{i\}$ and $H$ has order type $\beta_i$.
\end{definition}

\begin{lemma}\label{LEMtukeyantichain}
Suppose $\kappa=\cf\kappa>\omega$ and $P$ is a directed set such that $[\kappa]^{<\omega}\leq_T P$.  Then $P$ contains a set of $\kappa$\nbd-many pairwise incomparable elements.
\end{lemma}
\begin{proof}
Let $Q$ be a well\nbd-founded, cofinal subset of $P$.  Then $P\equiv_T Q$; let $f\colon[\kappa]^{<\omega}\leq_T Q$.  Define $g\colon[\kappa]^2\rightarrow 3$ by $g(\{\alpha<\beta\})=0$ if $f(\{\alpha\})\not\leq f(\{\beta\})\not\leq f(\{\alpha\})$ and $g(\{\alpha<\beta\})=1$ if $f(\{\alpha\})> f(\{\beta\})$ and $g(\{\alpha<\beta\})=2$ if $f(\{\alpha\})\leq f(\{\beta\})$.  By the Erd\"os\nbd-Dushnik\nbd-Miller Theorem, $\kappa\rightarrow(\kappa,\omega+1,\omega+1)$.  Since $Q$ is well\nbd-founded, there is no $H\in[\kappa]^\omega$ such that $g[[H]^2]=\{1\}$.  Since $f$ is Tukey and all infinite subsets of $[\kappa]^{<\omega}$ are unbounded, there is no $H\subseteq\kappa$ of order type $\omega+1$ such that $g[[H]^2]=\{2\}$.  Hence, there exists $H\in[\kappa]^\kappa$ such that $g[[H]^2]=\{0\}$; whence, $f[[H]^1]$ is a $\kappa$\nbd-sized, pairwise incomparable subset of $P$.
\end{proof}

\begin{theorem}
Suppose $\kappa=\cf\kappa>\omega$ and $X$ is a compactum such that every point has a local base not containing a set of $\kappa$\nbd-many pairwise incomparable elements.  Then some point in $X$ has $\pi$\nbd-character less than $\kappa$.
\end{theorem}
\begin{proof}
Combine Theorem~\ref{THMpichartukey} and Lemma~\ref{LEMtukeyantichain} to prove the contrapositive of the theorem.
\end{proof}

\begin{corollary}
Suppose $X$ is a compactum such that every point has a local base that is well quasi\nbd-ordered with respect to $\supseteq$.  Then some point in $X$ has countable $\pi$\nbd-character.
\end{corollary}

Finally, let us present a few results about local Noetherian type and topological embeddings.

\begin{lemma}\label{LEMocharcharsubspace}
Suppose $X$ is a space, $Y\subseteq X$, and $p\in Y$ satisfies $\character{p,Y}=\character{p,X}$.  Then $\ochar{p,X}\leq\ochar{p,Y}$.
\end{lemma}
\begin{proof}
Set $\lambda=\character{p,Y}$ and $\kappa=\ochar{p,Y}$; \wma\ $\lambda>\omega$ by Theorem~\ref{THMcharpibounds}.  By Lemma~\ref{LEMocharconverse}, we may choose a local $\la\lambda,\kappa\ra$\nbd-splitter $\mcA$ at $p$ in $Y$.  For each $U\in\mcA$, choose an open subset $f(U)$ of $X$ such that $f(U)\cap Y=U$.  Set $\mcB=f[\mcA]$.  Then $\card{\mcB}=\lambda$ because $f$ is bijective.  Suppose $\mcC\in[\mcB]^\kappa$ and $p$ is in the interior of $\bigcap\mcC$ with respect to $X$.  Then $p$ is in the interior of $Y\cap\bigcap\mcC$ with respect to $Y$, in contradiction how we chose $\mcA$.  Thus, $\mcB$ is a local $\la\lambda,\kappa\ra$\nbd-splitter at $p$ in $X$.  By Lemma~\ref{LEMocharconverse}, $\ochar{p,X}\leq\kappa$.
\end{proof}

\begin{definition}
For all infinite cardinals $\kappa$, let $u(\kappa)$ denote the space of uniform ultrafilters on $\kappa$.
\end{definition}

\begin{theorem}\label{THMocharomegainukappa}
For each $\kappa\geq\omega$, there exists $p\in u(\kappa)$ such that $\ochar{p,u(\kappa)}=\omega$ and $\character{p,u(\kappa)}=2^\kappa$.
\end{theorem}
\begin{proof}
Let $A$ be an independent family of subsets of $\kappa$ of size $2^\kappa$.  Set $B=\bigcup_{F\in[A]^\omega}\{x\subseteq\kappa:\forall y\in F\ \ \card{x\setminus y}<\kappa\}$.  Since $A$ is independent, we may extend $A$ to an ultrafilter $p$ on $\kappa$ such that $p\cap B=\emptyset$.  For each $x\subseteq\kappa$, set $x^*=\{q\in u(\kappa):x\in q\}$.  Then  $\{x^*:x\in A\}$ is a local $\la 2^\kappa,\omega\ra$\nbd-splitter at $p$.  Since $\character{p,u(\kappa)}\leq 2^\kappa$, it follows from Lemma~\ref{LEMocharconverse} that $\ochar{p,u(\kappa)}=\omega$ and $\character{p,u(\kappa)}=2^\kappa$.
\end{proof}

\begin{theorem}
Suppose $\kappa\geq\omega$ and $X$ is a space such that $\character{X}=2^\kappa$ and $u(\kappa)$ embeds in $X$.  Then there is an \omegaop\ local base at some point in $X$.  Hence, $\ochar{X}=\omega$ if $X$ is homogeneous.
\end{theorem}
\begin{proof}
Let $j$ embed $u(\kappa)$ into $X$.  By Theorem~\ref{THMocharomegainukappa}, there exists $p\in u(\kappa)$ such that $\ochar{p,u(\kappa)}=\omega$ and $\character{p,u(\kappa)}=2^\kappa$.  By Lemma~\ref{LEMocharcharsubspace}, $\ochar{j(p),X}=\omega$.
\end{proof}

\end{document}